\documentclass[proc]{edpsmath}
\usepackage{type1cm}        
\usepackage{comment}
\usepackage{graphicx}        
\usepackage{svg}
\usepackage{multicol}        
\usepackage[bottom]{footmisc}
\usepackage{placeins}
\usepackage{subcaption}
\usepackage{newtxtext}       %
\usepackage{xcolor}
\usepackage{amsthm}       %
\usepackage[varvw]{newtxmath}       
\usepackage{stmaryrd}	

\renewcommand{\textcolor}[2]{#2} 

\usepackage{soul}

\usepackage[overload]{empheq}

\usepackage{comment}

\def\x{{\bf x}}

\def\v{{\bf v}}

\def\bphi{{\boldsymbol \Phi}}

\def \h01{ \mathcal{H}_0^1}
\def\fflux{\mathbf{F}_\epsilon}

\def\dtn{{\rm DtN}}
\def\wdtn{ \widetilde{\rm DtN} }

\def\mass{a}

\def\dsp{\displaystyle}

\def\uoi{u_{ \O_i}}
\def\upoi{u_{\partial \O_i}}

\def\red{\textcolor{red}}
\def\mb{\textcolor{violet}}
\def\kb{\textcolor{magenta}}

\def\O{\Omega}
\def\G{\Gamma}

\def\RR{\mathbf{R}}

\def\node{s}
\def\triangles{{T}}

\def\Node{\mathfrak{s}}
\def\map{\mathbf{R}}
\def\mapH{\mathscr{R}}
\def\NH{\mathscr{N}}
\def\NHi{\NH_{\partial\Oi}}

\def\ov{\overline}

\def\Oc{{\ov{\O}}}

\def\Oic{{\ov{\O}_i}}
\def\Oi{{\O_i}}

\usepackage{multicol}
\usepackage{multirow}
\def\H10{H^1_{ \partial \O \setminus \partial \O_S}}

\def\dsp{\displaystyle}
\def\DtN{{\rm DtN}}

\def\fflux{\mathbf{F}}

\begin{document}

\title{Learning local Dirichlet-to-Neumann maps of nonlinear elliptic PDEs with rough coefficients
}\thanks{The authors are presented in alphabetical order.}
%

\author{Miranda Boutilier}\address{miranda.boutilier@univ-cotedazur.fr}
\author{Konstantin Brenner}\address{konstantin.brenner@univ-cotedazur.fr}
\author{Larissa Miguez}\address{lamiguez@posgrad.lncc.br}
%
%
\begin{abstract} 
Partial differential equations (PDEs) involving high contrast and oscillating coefficients are common in scientific and industrial applications.
Numerical approximation of these PDEs is a challenging task that can be
addressed, for example, by multi-scale finite element analysis. For
linear problems, multi-scale finite element method (MsFEM) is well
established and some viable extensions to non-linear PDEs are known.
However, some features of the method seem to be intrinsically based on
linearity-based. In particular, traditional MsFEM rely on the reuse of
computations. For example, the stiffness matrix can be calculated just
once, while being used for several right-hand sides, or as part of a
multi-level iterative algorithm. Roughly speaking, the offline phase of
the method amounts to pre-assembling the local linear
Dirichlet-to-Neumann (DtN) operators. We present some preliminary
results concerning the combination of MsFEM with machine learning tools.
The extension of MsFEM to nonlinear problems is achieved by means of
learning local nonlinear DtN maps. The resulting
learning-based multi-scale method is tested on a set of model nonlinear PDEs involving the $p-$Laplacian and degenerate nonlinear diffusion.
\end{abstract}

\maketitle

\section*{Introduction}

A variety of real-world phenomena are described by Partial Differential Equations (PDEs). When solutions to PDE models present a high variability in small spatial regions or short periods of time, this behavior is considered  multiscale. 
Standard numerical methods such as Finite Element, Finite Difference, or Finite Volume methods often have difficulty approximating multiscale behavior.
Therefore, multiscale numerical methods have emerged as attractive options for dealing with such problems, including  the Localized Orthogonal Decomposition  method \cite{lod}, Heterogeneous Multiscale methods \cite{hmm}, the Generalized Finite Element Method \cite{gfem}, the Multiscale Spectral Generalized Finite Element Method\cite{msgfem, msgfem2} and the Multiscale Finite Element Method (MsFEM) \cite{hou1997multiscale}. MsFEM methods appear as a two-scale finite element approach: a macroscopic scale to capture the global solution behavior of the material and a microscopic scale to model fine heterogeneities. By integrating information from these two scales, MsFEM methods provide a more accurate solution for  problems with multiscale behavior.

Furthermore, it is worth noting these methods offer the advantage of reusing calculations to enhance efficiency. For instance, the stiffness matrix of MsFEM can be computed once and reused multiple times, just as Domain Decomposition (DD) \cite{ddbook} techniques can benefit from the precomputation of LU decompositions of local matrices. However, MsFEM strategies do not naturally extend to nonlinear problems ~\cite{chaturantabut2010nonlinear, efendiev2013generalized}; we focus therefore  on the the computation of local Dirichlet-to-Neumann (DtN) operators. 
Specifically, with coarse (macroscopic) and fine (microscopic) partitionings of the heterogeneous domain, we can compute local DtN operators on each subdomain. These DtN local operators, to be defined in detail later, provide a relationship between Dirichlet and Neumann boundary data; these local operators lead to a nonlinear substructuring method such that the fine-scale discrete
 nonlinear problem can be written in terms of the local DtN maps.

On the other hand, with the renewed interest in Machine Learning (ML) within the scientific computing community, numerous techniques have been proposed to apply ML techniques to the  solution of PDEs. In \cite{lu2021learning}, the authors employ DeepONet neural networks to learn nonlinear operators. Additionally, the authors of \cite{badia2023finite} introduce an approach that combines Finite Element interpolation techniques with Neural Networks (NN) to solve a wide range of direct and inverse problems. Among these approaches are Physics-Informed Neural Networks (PINNs) \cite{raissi2019physics}. The central idea behind PINNs is to minimize a functional representing the residual of the PDE and its initial and boundary conditions. Since then, many pieces of work have been published on PINNs, including \cite{mao2020physics, cai2021physics, misyris2020physics}. However, problems with complex domains have led to other methodologies based on DD methods and PINNs, including Extended PINNs \cite{jagtap2021extended}, Conservative PINNs \cite{jagtap2020conservative}, Variational PINNs \cite{kharazmi2019variational} and Finite Basis Physics-Informed Neural Networks \cite{moseley2023finite}. Additionally,  \cite{huge2020differential} proposed Differential Machine Learning (DML),  originally with a financial application. 
In DML, supervised learning is ``extended" such that ML models are trained not only with
values/inputs but also and on correpsonding derivatives. 

In this work,  we aim to use machine learning techniques to learn the local nonlinear DtN maps. The neural network is trained to replicate the action of the local nonlinear DtN maps on some coarse subset of the trace space. Once the training is completed, the surrogate DtN operators will be used to solve a nonlinear substructuring method, which is solved via Newton's method.
For training the nonlinear maps, we use an approach similar to DML; that is, we  incorporate gradient information into the loss function. Additionally, we impose a monotonicity property into the loss function which is based on the known monotonicity of the DtN maps. 

This paper is organized as follows. In Section \ref{sec:dtn}, we introduce the model problem and associated approximate substructured formulation, which includes the definition of the continuous DtN maps. In Section \ref{sec:dtncomp}, we introduce the finite element discretization, discrete substructured formulation, and discrete DtN maps. We also provide a description of the machine learning process which is used to determine the learned DtN operators and assemble the learned substructured formulation. Section \ref{sec:numericalExp} contains numerical experiments, showing both the accuracy of the learned DtN maps and the resulting learned approximate substructured formulation.  Section \ref{sec:conclusion} concludes with a summary and discussion for future work.

\section{Model Problem and Approximate Substructured Problem}\label{sec:dtn} 
Let $\O$ be either an open connected polygonal domain in $ \mathbb{R}^2$ or a real interval. In this work, we focus on model problems of the form
	\begin{equation}\label{model_pde}
		\left\{
		\begin{array}{rll}
			\mass u  + \text{div} \left( K_\epsilon(\mathbf{x}) \fflux( u, \nabla u) \right) &=& 0  \qquad \mbox{in} \qquad \O, \\
			u &=& u_D \qquad \mbox{on} \qquad  \partial \O,\\
		\end{array}
		\right.
	\end{equation} 
where $\mass$ is a positive real coefficient and  $K_\epsilon(\mathbf{x})$, representing the heterogeneity of the problem, is assumed to satisfy $K_\epsilon(\mathbf{x}) \geq k_{min} >0$ for all $\x \in \O$. For the ``flux functions'' $F(u, \nabla u)$, we consider either the model
$$
F(u, \nabla u) = \kb{-}\nabla u^p, \quad p > 1,
$$
leading to a porous media problem, or the $p-$Laplace model
$$
F(u, \nabla u) = \kb{-}|\nabla u|^p \nabla u, \quad p > 0.
$$
The well posedness of the problem \eqref{model_pde} is known (see e.g. 
\cite{vazquez2007porous} and \cite{leray1964quelques}).

\subsection{Substructured and Approximate Substructured Problem}
\begin{figure}
    \centering
    \includegraphics[height=5cm]{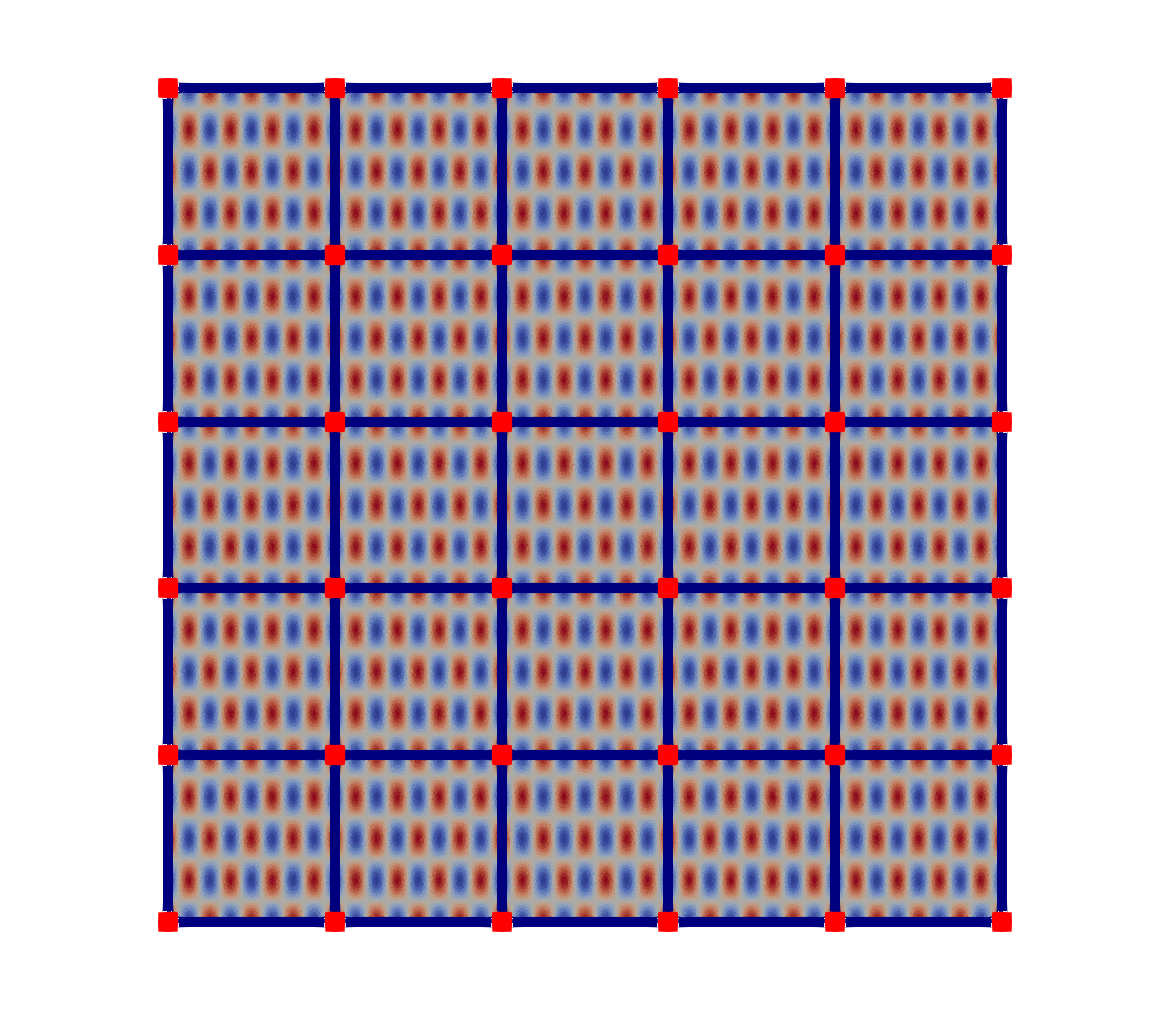}
    \caption{Heterogeneous domain $K_\epsilon(\mathbf{x})$ with $N=5 \times 5$ coarse cells. Nonoverlapping skeleton $\Gamma$ is denoted by thick dark blue lines.}
    \label{fig:dofs}
\end{figure}

Consider a finite nonoverlapping rectangular partitioning of $\O$ denoted by $\left( \O_j \right)_{j= 1}^{N}$. We will refer to  $\left( \O_j \right)_{j= 1}^N$  as the coarse mesh over $\O$. Additionally, we denote by $\Gamma$ its skeleton,  that is  $\Gamma=  \bigcup_{j = 1}^N \partial \O_j$. With this, each edge of the skeleton $\Gamma$ is denoted by $\Gamma_{ij}= \partial \O_i \cap \partial \O_j$,
while the set of the coarse grid nodes  is defined by $\{ \Node_\alpha \}_{\alpha=1}^{\NH_\G} = \bigcup_{i,j = 1}^{N} \partial \G_{ij}$.
An example of the coarse partitioning and corresponding coarse grid nodes is shown in Figure \ref{fig:dofs}.

The Dirichlet-to-Neumann (DtN) map provides a mathematical relationship between Dirichlet and Neumann boundary data,  allowing us to compute one set of boundary data from the other, which can be very useful in solving boundary value problems. Given a nonoverlapping partitioning of domain, these DtN maps can be defined and computed locally on each subdomain. \kb{Let us consider a local Dirichlet problem 
	\begin{equation}\label{eq:localprobs}
	\left\{
	\begin{array}{rll}
			 \mass u_i  + \text{div} \left( K_\epsilon(\mathbf{x}) \fflux( u_i, \nabla u_i) \right) &=& 0  \qquad \mbox{in} \qquad \O_i, \\
		\dsp  u_i &=& g_i \qquad \mbox{on} \quad   \partial \O_i.
	\end{array}
	\right.
\end{equation}
For a given $g_i$ we formally denote by $\dtn_{i}(g_i)$ the  the Neumann traces of $u_i$ on $\partial \Oi$, that is
$$
\dtn_{i}(g_i) = K_\epsilon(\mathbf{x})\fflux(u_i, \nabla u_i)|_{\partial \O_i} \cdot \mathbf{n}_i \qquad \text{on} \quad \partial\Oi,
$$
where $\mathbf{n}_i$ denotes the unit outward normal on $\partial \Oi$.}

Our main idea with the nonlinear substructuring method is to replace the original problem \eqref{model_pde} by the one which only involves the unknown Dirichlet data $g$ defined over $\G$. This new unknown function is determined by imposing a matching flux condition
along the edges of the skeleton $\Gamma$. 




For \kb{a regular enough function $g$ defined on $\G$}, we formally introduce the local jump operator, returning the mismatch in fluxes on a given coarse edge
\kb{ $$\llbracket \rrbracket_{ij}: g \mapsto  \dtn_{i}(g|_{\partial \O_i})|_{\Gamma_{ij}} + \dtn_{i}(g|_{\partial \O_j})|_{\Gamma_{ij}}.$$
}
and the global jump operator $\llbracket \rrbracket$ is given by 
$ \llbracket g \rrbracket |_{\G_{ij}}=  \llbracket g\rrbracket_{ij}$.
Enforcing zero flux jump at every coarse edge $\G_{ij}$, leads to a global substructured problem that can be written in a weak form as: find $g$ satisfying $g|_{\partial \O} = u_D$ \kb{such that for all regular enough $w$ vanishing on $\partial \O$ it holds that
\begin{align} \label{eq:finesubstruct}
	\int_\G \llbracket g \rrbracket w \, d\gamma(\x) = 0. 
	\end{align} 
}
Let us assume for the moment 
that \kb{we are in the possession of a procedure allowing to compute evaluate the} DtN operators. Then, the problem \eqref{eq:finesubstruct} can be discretized by \kb{taking the test and trial function 
from a}
finite dimensional space. 
\kb{In particular, we consider the space $V_H(\G)$ made of the continuous functions on $\G$ that are  affine on every coarse edge of $\G$, that is
\begin{equation}\label{eq:VH}
V_{H}(\Gamma)= \{ v \in C^0(\Gamma) \quad \text{s.t.} \quad v|_{\Gamma_{ij}} \in \mathbb{P}_1(\Gamma_{ij})  \quad \text{ for all } \quad  \G_{ij} \}.
\end{equation}
}
Let $V_{H,0}(\Gamma)= \{ v\in V_{H}(\Gamma) \quad \text{s.t.} \quad v|_{\partial{\O}} = 0\}$ and 
let $u_{D,H} \in V_H(\G)$ denote some approximation of $u_D$, e.g. obtained by interpolation. To form the approximate substructured problem in $V_H(\G)$, one then replaces the problem \eqref{eq:finesubstruct} by: find $g_H$ satisfying $g_H|_{\partial \O} = u_{H,D}$ and 
\begin{align} \label{eq:coarsesubstruct}
	\int_\G \llbracket g_H \rrbracket w_H \, d\gamma(\x) = 0 \qquad \forall w_H \in V_{H,0}(\G). 
	\end{align} 


\def\IS{{\cal N}}
\section{Discrete Finite Element Substructuring and DtN Maps}\label{sec:dtncomp}
\subsection{Finite element discretization and domain decomposition}\label{sec:fem}
Let $\triangles_h$ denote the triangulation of $\O$ which is supposed to be conforming with respect to partitioning $\left( \O_j \right)_{j= 1}^{N}$. Let 
\mb{$\{\node_k\}_{k=1}^{\ov \O}$} be the set of vertices (or points) of $\triangles_h$. We introduce the finite element space defined by
$$
V_h =  \{ v  \,\,|  \,\, v \in C^0(\overline{\O}), \quad \text{s.t.} \quad v|_{t} \in \mathbb{P}_1, \quad \text{for all} \quad  t \in \triangles_h \},
$$
and we denote the corresponding nodal ``hat" basis functions $(\eta^{k})_{k=1}^{N_\O}$ for each point in $\node_k$. Let $V_{h,0}= V_h \cap H_0^1$, and
let $u_{h,D}\in V_h$ be some approximation of the boundary data $u_D$.
The Galerkin finite element formulation is as follows: 
\begin{equation}\label{eq:galerkin}
\begin{array}{lllll}
\text{Find $u_h \in V_{h}$ satisfying $u_h|_\O= u_{h,D}$ such that} \\[2pt]
\dsp \int_{\Omega} \,  \mass u_h v_h \kb{-} (K_{\epsilon}(\mathbf{x}) \fflux(u_h, \nabla u_h)) \cdot \nabla v_h  \, d\x = 0 \qquad \text{for all} \quad v_h \in V_{h,0}.
\end{array}
\end{equation}
Now, consider the following linear reconstruction operator 
$\pi_h$ acting from $\mathbb{R}^{N_\O}$ to $V_h$, such that 
$$
\dsp \pi_h(v)(\x)= \sum_{k=1}^{N_\O} v^{k} \eta^k (\x). 
$$
We note that $\pi_h$ is bijective, with $\left(\pi^{-1}(\v_h)\right)^k = v_h(s_k)$ for any triangulation point $s_k$.

We denote by $F$ the ``Neumann residual'' function corresponding to \eqref{eq:galerkin}; more precisely $F$ acts from $\mathbb{R}^{N_\O}$ to itself and has its $l$th component defined by 
\begin{equation}\label{eq:disc}
\left( F(u) \right)^l=  \int_{\O} \pi
_h(u) \, \eta^l \kb{-} K_\epsilon(\mathbf{x})  \fflux( \pi_h(u), \nabla \pi_h(u)) \cdot \nabla \eta^l  \,d\x,
\end{equation}
for $l = 1,\ldots, N_\O$.

Before moving any further, let us introduce some notation related to the decomposition of the mesh and finite element degrees of freedom induced by the partitioning $\left(\Oi\right)^{N}_{i =1}$ of $\O$. Let  $J = \{1, \ldots, N_{\ov \O} \}$ denote the index set of the triangulation points, for any $\omega \in \ov \O$ we introduce a subset of indices $$J_\omega = \left\{ j_\omega^1, \ldots j_\omega^{N_\omega}\right\} \subset J$$ such that 
    $j\in J_\omega$ if and only if the point $\node_j$ belongs to $\omega$. We note that the number of elements in $J_\omega$ is denoted by $N_\omega$.
    
The standard $N_\omega$ by $N_{\ov\O}$ boolean restriction matrix $\RR_\omega$ ``from $\ov \O$ to $\omega$'' is defined by
    $$
    \left( \RR_\omega \right)_{kl} = \left\{
    \begin{array}{lllll}
       1  & \text{if} \quad l = j_\omega^k,  \\
       0  & \text{else,}
    \end{array}
    \right.
    $$
for any $k = 1,\ldots, N_\omega$ and $l = 1, \ldots, N_{\ov \O}$. Further, for any $A, B \subset \ov \O$, we denote $\map^A_B = \RR_A (\RR_B)^T$. If $A\subset B$ the matrix $\map^A_B$ can be interpreted as a restriction matrix ``from $B$ to $A$''. Conversely, if $B\subset A$, the same $\map^A_B$ may be interpreted as an extension matrix ``from $B$ to $A$''.  Let us mention a couple of useful properties of this restriction/extension linear operator. First, we note that for any $C$ such that $A, B\subset C \subset \ov \O$ we have 
$$
\map^A_C \map^C_B = \map^A_B = \map^A_{A \cap B}\map^{A\cap B}_B.
$$

With this, the discrete problem \eqref{eq:galerkin} can be expressed as follows: 
\begin{equation}\label{eq:pbalg}
\begin{array}{llll}
\text{Find $u \in \mathbb{R}^{N_\O}$ such that $\map_{\partial \O} u = \kb{\map_{\partial \O}}\pi_h^{-1}\left( u_{h,D} \right) $ and satisfying} \,\,
\map_{\O} F(u) = 0.
\end{array}
\end{equation}

We note that the integral in \eqref{eq:galerkin} can be split into a sum of integrals over $\Oi$, meaning that the finite element residual can be assembled by gluing together some  the local contributions. We detail below this local assembly procedure because of it use for the substructured formulation. 

For a subdomain $\Oi$ let $\{\eta_i^l\}_{l = 1 }^{N_{\Oic}}$ denote the set of local finite element basis functions such that $\eta_i^l(\x) = \eta^{j_\Oic^l}(\x)|_\Oi$ for all $l = 1,\ldots, N_{\Oic}$. The local function reconstruction operator $\pi_{h,i}$ is defined by $\dsp \pi_{h_i}(v_i)(\x) = \sum_{l = 1}^{N_{\Oic}} v^l \eta_i^l(\x)$. Similar to \eqref{eq:disc} we define the local residual functions $F_i$. That is for any $u_i \in \mathbb{R}^{N_\Oic}$ and  $l=1, \ldots, N_\Oic$, 
we set
\begin{equation}\label{eq:localdisc}
\left( F_i(u_i) \right)^l=  \int_{\O_i} \pi
_{h,i}(u) \, \eta^l_i \kb{-} K_\epsilon(\mathbf{x})  \fflux( \pi_{h,i}(u_i), \nabla \pi_{h,i}(u_i)) \cdot \nabla \eta^l_i  \,d\x. 
\end{equation}
with $\eta_i^l = \eta^{j_\Oic^l}$. For any $u \in \mathbb{R}^{N_\Oc}$ the global residual $F(u)$ is obtained from the local components by the following expression
\begin{equation}\label{eq:neumannsum}
F(u) = \sum_{i = 1}^N \RR_{\Oic}^T F_i(\RR_\Oic u ).
\end{equation}
.

\subsection{Discrete $\dtn$ maps}
In this section  we introduce the finite element version the Dirichlet-to-Neumann operator, associated with our model PDE and a subdomain $\Oi$. This discrete operator, denoted by $\dtn_{h,i}$, is defined as the nonlinear mapping from the set of degrees of freedom  $\mathbb{R}^{N_{\partial \O_i}}$ associated with the boundary of $\Oi$ to itself. We also provides the formula for the (Fr{\'e}chet) derivative of $\dtn_{h,i}$. The construction of $\dtn_{h,i}$ relies on the definition \eqref{eq:localdisc} of the local ``Neumann'' residual $F_i$.

Let us split \mb{the set of induces} $J_{\Oic}$ into two non-overlapping subsets $J_\Oi$ and $J_{\partial \Oi}$ associated with internal degrees of freedom located in $\O_i$ and on the boundary $\partial \O_i$. 
This gives the representation of the vector $u_i \in \mathbb{R}^{N_{\Oic}}$ associated to the subdomain $\Oic$ in the form $u_i= \map_{\Oi}^{\Oic} \uoi + \map_{\partial \Oi}^{\Oic} \upoi$ with some $u_\Oi \in \mathbb{R}^{N_\Oi}$ and $u_{\partial \Oi} \in \mathbb{R}^{N_{\partial \Oi}}$.
In turn, we express the local Neumann residual as following
\begin{equation}\label{eq:Fi_split}
F_i( u_i ) = 
\map_{\Oi}^\Oic F_\Oi (\uoi, \upoi) + 
\map_{\partial \Oi}^\Oic F_{\partial \Oi} (\uoi, \upoi)
\end{equation}
where 
$$
F_{\O_i}(\uoi, \upoi) = \map_{\Oic}^{\Oi} F_i\left( \map_{\Oi}^{\Oic} \uoi + \map_{\partial \Oi}^{\Oic} \upoi \right)
$$
and
$$
F_{\partial \O_i}(\uoi, \upoi) = \map_{\Oic}^{\partial\Oi} F_i\left( \map_{\Oi}^{\Oic} \uoi + \map_{\partial \Oi}^{\Oic} \upoi \right).
$$
The discrete counterpart of the local Dirichlet problem \eqref{eq:localprobs} can be expressed as: 
\begin{equation}\label{eq:localprobs_fem}
\text{Find $\uoi$ such that} \,\,
F_{\O_i}(\uoi, g_{\partial \Oi}) = 0
\end{equation}
for a given $g_{\partial \Oi} \in \mathbb{R}^{N_{\partial \Oi}}$.
Assuming that \eqref{eq:localprobs_fem} admits unique solution for any $g_{\partial \Oi}$, we define the solution operator $G_{\Oi}$ from $\mathbb{R}^{N_{\partial \Oi}}$ to $\mathbb{R}^{N_{\Oi}}$ such that 
\begin{equation}\label{eq:solop}
F_{\Oi}\left( G_{\Oi}( g_{\partial \Oi} ),  g_{\partial \Oi}\right) = 0
\end{equation}
for all $g_{\partial \Oi}$. With that, the discrete DtN operator is defined by
\begin{equation}\label{eq:DtNhi}
\dtn_{h,i}(  g_{\partial \Oi} ) = F_{\partial \Oi}\left( G_{\Oi}( g_{\partial \Oi} ),  g_{\partial \Oi}\right).
\end{equation}

Under the assumptions of  Implicit Function Theorem (see e.g. Proposition 5.2.4. of \cite{OR70}) the Fr{\'e}chet derivative of $\dtn_{h,i}$ can be expressed as
\begin{align}\label{eq:dtnprimeold}
\dtn_{h,i}'(\upoi) & = \partial_{\O_i} F_{\partial \O_i}( G_{\O_i}(\upoi), \upoi)G'_{\O_i}(\upoi)      + \partial_{ \partial \O_i} F_{ \partial \O_i}(G_{\O_i}(\upoi), \upoi). 
\end{align}
Recalling that $u_{\O_i}= G_{\O_i}(\upoi)$, we have 
\begin{align}\label{eq:dtnprimegood}
\dtn_{h,i}'(\upoi) & = \partial_{\O_i} F_{ \partial \O_i}( \uoi,\upoi)G'_{\O_i}(\upoi)      + \partial_{ \partial \O_i} F_{ \partial \O_i}(\uoi, \upoi), 
\end{align}
where, thanks to the identity 
\begin{equation}\label{eq:gdef}
    F_{\O_i}(G_{\O_i}(\upoi), \upoi)=0,
\end{equation}  
we have 
  \begin{align}\label{eq:gprime}
 G'_{\O_i}(\upoi)  &= -\left(  \partial_{\O_i}  F_{ \O_i}(\uoi, \upoi) \right)^{-1} \partial_{ \partial \O_i}  \left(F_{\O_i}(\uoi,\upoi)\right). 
 \end{align}
Plugging \eqref{eq:gprime} into \eqref{eq:dtnprimeold}  results in
\begin{align}
\dtn_{h,i}'(\upoi) & = -\partial_{\O_i} F_{\partial \O_i}(\uoi, \upoi) \left(  \partial_{\O_i}  F_{ \O_i}(\uoi, \upoi) \right)^{-1} \partial_{ \partial \O_i}  F_{ \O_i}(\uoi, \upoi)  
+ \partial_{ \partial \O_i} F_{ \partial \O_i}(\uoi, \upoi), \nonumber
\end{align}
We remark that computing $u_{\O_i}$ requires solving the local nonlinear problem \eqref{eq:localprobs_fem} with $\upoi$ as data. This is typically done my some fixed point method. Ones the approximate value of $u_{\O_i}$ is obtained, computing $\dtn_{h,i}'(\upoi)$ requires solving the linear system with matrix right-hand-sides given by  
$F_{ \O_i}(\uoi, \upoi)$. If $u_{\O_i}$ is obtained by exact Newton's method, the former can be done at marginal computational cost. 


\subsection{Discrete Substructured Problems}\label{sec:disc }
With $\dtn_{h,i}$ introduced above, let us provide the substructured formulation of the discrete problem \eqref{eq:pbalg}. 
Let 
$g\in\mathbb{R}^{N_\G}$ be a vector representing the unknown solution values at the skeleton $\G$. Denoting $g_{\partial \Oi} = \map^{\partial \Oi}_\G g$ and $u_i = G_\Oi\left( g_{\partial \Oi}  \right)$,
we express the unknown discrete solution as
\begin{equation}\label{ansatz}
u = \RR_\G^T g + \sum_{i = 1}^N \RR_\Oi^T u_\Oi.
\end{equation}
Observing that $\RR_{\Oic} u = \map^\Oic_{\partial \Oi}g_{\partial \Oi} + \map^\Oic_\Oi \kb{u_\Oi}$,  in view of \eqref{eq:Fi_split},  we deduce that
\begin{equation}\label{some_eq_1}
F_i\left(  \RR_{\Oic} u  \right) = 
\map_{\Oi}^\Oic F_\Oi (u_\Oi, g_{\partial \Oi}) + 
\map_{\partial \Oi}^\Oic F_{\partial \Oi} (u_\Oi, g_{\partial \Oi}).
\end{equation}
The first term in the right hand side of \eqref{some_eq_1} is zero by definition of $u_\Oi$, which yields
$$
F(u) = \sum_{i = 1}^N \RR_{\partial \Oi}^T   \dtn_{h,i}\left( g_{\partial \Oi} \right) = \sum_{i = 1}^N \RR_{\partial \Oi}^T   \dtn_{h,i}\left( \map^{\partial \Oi}_\G g\right), 
$$
in view of \eqref{eq:neumannsum} and \eqref{eq:DtNhi}. Let us denote
\begin{equation}\label{eq:FG_splitting}
    F_\G(g) = \sum_{i = 1}^N \map^\G_{\partial \Oi}   \dtn_{h,i}\left( \map^{\partial \Oi}_\G g\right).
\end{equation}
Substituting $F(u) = \RR_\G^T F_\G(g)$ in \eqref{eq:pbalg} leads to the equation 
\mb{$\map_\O \map_\G F_\G(g) = 0$.}
Since the latter is trivially satisfied for all degrees of freedom not lying on $\G$, the system has to be reduced to $\G \cap \O$, leading to
\mb{$$
0 = \map^{\G\cap\O}_\O \map_\O \map_\G  F_\G(g) = \map^{\G\cap\O}_\G  F_\G(g).
$$}
$$$$
We have obtain the following substructured problem: find $g\in \mathbb{R}^{N_\G}$ such that $\map^{\partial \O}_\G g = \RR_{\partial \O}\pi_h^{-1}(u_{h,D})$  and 
\begin{equation}\label{eq:nonlinfinesubstruct}
\map^{\G\cap\O}_\O F_\G(g) = 0.
\end{equation}

\subsection{Coarse discrete $\dtn$ maps and discrete approximate substructuring}
We now proceed to the discrete version of the  approximate substructured problem \eqref{eq:coarsesubstruct}. 
Similar to the definitions introduced in Section \ref{sec:fem}, for any $\omega \subset \G$, we denote by $\NH_{\omega}$ the number of coarse nodes belonging to $\omega$, and by $\mapH^\omega_\G$ the restriction matrix from the
$\mathbb{R}^{\NH_\G}$ to $\mathbb{R}^{\NH_\omega}$. As before, the transpose of $\mapH^\omega_\G$ is denoted by $\mapH_\omega^\G$. 

Let 
$(\phi^\alpha)_{ \alpha=1}^{\NH_\G}$ be a nodal  basis of $V_H(\Gamma)$ made of ``hat'' functions associated to the set coarse grid nodes $\{ \Node_\alpha \}_{\alpha=1}^{\NH_\G}$ (see Figure \ref{fig:dofs} for visual). More precisely, the basis function $\phi^\alpha \in V_H(\G)$ is assumed to satisfy
	\[ \phi^\alpha(\Node_j)=\begin{cases} 
		1, & \alpha=j,   \\
		0, & \alpha \neq j, \\
	\end{cases}
	\]
 for all coarse grid nodes $\Node_j$. Given this nodal basis in $V_H(\G)$ we define the function reconstruction operator $\pi_H$ from $\mathbb{R}^{\NH_\G}$ to $V_H(\G)$ by $\dsp \pi_H(v)(\x) = \sum_{\alpha=1}^{\NH_\G} v^\alpha \phi^\alpha(\x)$.
We further denote by $\bphi_H$ the $N_\Gamma$  by $\NH_\G$ \mb{matrix} whose columns contains the values of the basis functions $\phi^\alpha$ at the skeletal degrees of freedom; more specifically, \mb{the $j$th component of a column index $\alpha$ of $\bphi_H$ is given by} 
$$
(\bphi_H)^{j,\alpha} = \phi^\alpha(\node_{j_\G^j}).
$$
for $\alpha=1, \ldots, \NH_\G$ and each $j=1, \ldots, N_\G$. 
%
%

The local version of $\bphi_{\mb{H}}$ is defined by
\begin{equation}
\bphi_{H i} = \map_{\Gamma}^{ \partial \O_i} \bphi_{H}
\mapH_{\partial \Oi}^\G,
\end{equation} 
where $\bphi_{Hi}$ is composed of the $\NHi$ local basis vectors of size 
$\mb{N_{\partial \O_i}}$ corresponding to the local coarse grid nodes on $\partial \O_i$. 
We note that $\bphi_{H i}$ satisfies
\begin{equation}\label{some_eq_2}
   \map_{\Gamma}^{ \partial \O_i} \bphi_{H} = \bphi_{Hi} \mapH^{\partial \Oi}_\G. 
\end{equation}
In order to derive the coarse approximation of the problem \eqref{eq:nonlinfinesubstruct} we first express the latter in the following variational form:
\begin{equation}\label{eq:substruct_fine_var}
\begin{array}{llll}
\text{Find $g$ such that $\map_\G^{\partial \O} g = \RR_{\partial \O}\pi_h^{-1}(u_{h,D})$ and satisfying}\\[7pt]
\dsp v^T \, F_\G(g) = 0 \qquad \text{for all} \quad v\in 
{\rm Im} \left( \map_{\G\cap \O}^\G \right).
\end{array}
\end{equation}
The Galerkin approximation of \eqref{eq:substruct_fine_var} using the coarse space $V_H(\G)$ is obtained by setting $g \approx \bphi_{H} g_H$ with some appropriate boundary conditions and taking the test elements $v$ in 
${\rm Im} \left( \map_{\G\cap \O}^\G \right) \cap {\rm Im\left( \bphi_{H} \right)} = {\rm Im\left( \bphi_{H} \right)} \cap {\rm Im} \left( \mapH_{\G\cap \O}^\G \right)$. The latter amount to set $v = \bphi_H \mapH_{\G\cap \O}^\G \v_H$ with $v_H\in \mathbb{R}^{\mathfrak{N}_{\G\cap\O}}$. This leads to the system 
$$
 \mapH^{\G\cap \O}_\G \bphi_H^T F_\G\left( \bphi_H g_H \right) = 0.
$$
In view of \eqref{eq:FG_splitting} and \eqref{some_eq_2}, we have
$$
 \sum_{i = 1}^N \bphi_H^T \map^\G_{\partial \Oi}   \dtn_{h,i}\left( \map^{\partial \Oi}_\G \bphi_H^T  g_H\right) 
 = 
 \sum_{i = 1}^N \mapH_{\partial \Oi}^\G  \bphi_{Hi}^T  \dtn_{h,i}\left( \bphi_{Hi} \mapH^{\partial \Oi}_\G g_H\right).
$$
With this, the definition of the coarse  $\dtn$ operator acting from $\mathbb{R}^{\NHi}$ to itself given by 
\begin{equation}\label{eq:dtnHi}
    \dtn_{H,i} = \bphi_{Hi}^T \circ \dtn_{h,i} \circ \bphi_{Hi},
\end{equation}
such that we denote
\begin{equation}\label{eq:FH}
F_H(g_H) =  \sum_{i = 1}^N \mapH_{\partial \Oi}^\G  \dtn_{H,i}( \mapH^{\partial \Oi}_\G g_H).
\end{equation}
Let $u_{H, D} \in V_H(\G)$ be some approximation of $u_D$, the coarse Galerkin approximation of \eqref{eq:pbalg}, or equivalently, the finite element approximation of \eqref{eq:coarsesubstruct} reads as follows: 
\begin{equation}\label{eq:nonlincoarsesubstruct}
\begin{array}{llll}
\text{Find $g_H \in \mathbb{R}^{\NH_\G}$ such that $\mapH_{\G}^{\partial \O} g_H = \mapH_{\G}^{\partial \O} \pi_H^{-1}(u_{H, D})$ and}\,\,  
\mapH_{\G}^{\G\cap\O} F_{H}(g_H) = 0.
\end{array}
\end{equation}


We remark that even though it is considered a coarse approximation, the computation of the discrete $\text{DtN}_{i}$ relies on fine-scale information and fine-scale input. That is, it still takes in as input a vector of size $N_\Gamma$ which lives on the fine degrees of freedom along the skeleton $\Gamma$. This computation can still become quite expensive and the cost is not to be negated, as $\text{DtN}_{i}$ must be computed at each Newton iteration used to solve the discrete substructured problem.

\subsection{Learning Dirichlet-to-Neumann maps}\label{sec:learning}

In Section \ref{sec:dtncomp}, we discussed the computation of the DtN map and how by nature, the substructured problem results in the DtN map being computed at each Newton iteration. As this computation can get very expensive, we explore a Scientific Machine Learning (SciML) application to our problem. Specifically, we introduce learned local DtN operators which we denote as $ \widetilde{\text{DtN}}_{H,i} $.Once the model is generated, we can call the model at each Newton iteration without the need for a separate computation on each subdomain. We remark that if the heterogeneities are periodic in each subdomain, we can use the same NN model.


Let $\widetilde{\text{DtN}}_{H,i}$ denote some approximation of $\text{DtN}_{H,i}$ to which we will refere to as to a surrogate model, and which is going to be introduced in more details below. By replacing in \eqref{eq:FH} 
the original operator $\text{DtN}_{H,i}$ by $\widetilde{\text{DtN}}_{H,i}$ we define
\kb{
$$
F_{H,learn}(g_H) =  \sum_{i = 1}^N \mapH_{\partial \Oi}^\G  \wdtn_{H,i}( \mapH^{\partial \Oi}_\G g_H).
$$
and we obtain the problem:
\begin{equation}\label{eq:approxsubstlearn}
\begin{array}{llll}
\text{Find $\widetilde{g}_H \in \mathbb{R}^{\NH_\G}$ such that $\mapH_{\G}^{\partial \O} \widetilde{g}_H = \mapH_{\G}^{\partial \O} \pi_H^{-1}(u_{H, D})$ and}\,\,  
\mapH_{\G}^{\G\cap\O} F_{H,learn}(\widetilde{g}_H) = 0
\end{array}
\end{equation}
}

\begin{figure}[ht!]
    \centering
    \includegraphics[width=.8\textwidth]{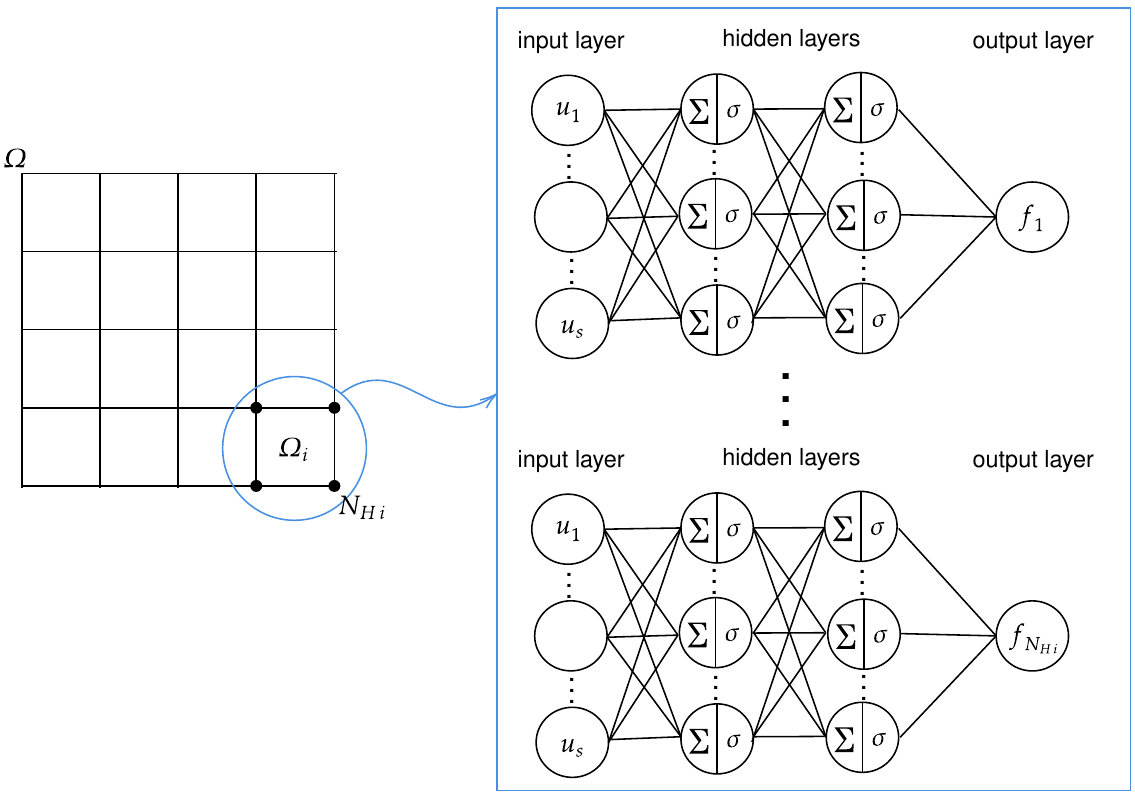}
    \caption{Representation of the learning workflow for the operator. Each neural network is depicted as a graph and will be responsible for learning each component of the $\widetilde{\text{DtN}}_{H,i}$ operator. The first column of neurons (from left to right in NN) is the input layer, taking inputs $ (u_1, u_2, \cdots, u_s)$ , and the last column is the output layer with output $f_{N_{H i}} \approx \widetilde{\text{DtN}}_{H,i}$. The intermediate columns represent the hidden layers.}
    \label{fig:network}
\end{figure}

We will construct the surrogate mapping $\widetilde{\text{DtN}}_{H,i}$ by learning the individual components of the original map. To learn this operator, we are using a fully-connected feed-forward neural network (NN). These networks are built as a sequence of layers including an input layer, hidden layers, and an output layer. Each layer performs an affine transformation followed by a nonlinear activation function $\sigma : \mathrm{R} \mapsto \mathrm{R}$. For each component $\widetilde{\text{DtN}}_{H,i}^l$, where $l = 1,\ldots, \NHi$, we approximate it using a neural network expressed as a linear combination of functions:
\begin{equation*}
    \widetilde{{\rm DtN}}_{H,i}^{l}(\theta; u) = \varphi_\theta^k \circ \ldots \circ \varphi_\theta^1(u),
\end{equation*}
where $\varphi_\theta^k(u)= \sigma (W_k(\theta) u + \mathbf{b}_k(\theta) )$. The weights matrix $W_k$ and biases $\mathbf{b}_k$ are determined through training. We illustrate this concept with Figure \ref{fig:network}.






In NN methodology, the solution to the DtN operator is approximated by the following optimization process. 
Let $ U= \{u^s\}_{s=1}^{n_s}$ be a set of sampling vectors sampled in the range $(u_{\text{min}}, u_{\text{max}})$, where $n_s$ denotes the total number of sampling vectors.
For each $i=1, \ldots, N$ and $l=1, \ldots \NHi$,  find $\theta^* = \underset{\theta}{\mathrm{arg\,min}}\,\mathcal{L}(\theta,U)$, with

    \begin{align}\label{eq:lossFunctional}
    \mathcal{L}_{i}^l(\theta, U) &= c_0 \mathcal{L}_{i,0}^l(\theta, U) + c_1 \mathcal{L}_{i,1}^l(\theta, U) + c_{\text{mon}} \mathcal{L}_{i,\text{mon}}^l(\theta)
    \end{align}
where $c_{0}$, $c_{1}$ and $c_{\text{mon}}$ are positive weights, $\mathcal{L}_{i,0}^l(\theta, u)$,  $\mathcal{L}_{i,1}^l(\theta, u)$ and $\mathcal{L}_{i,\text{mon}}^l(\theta)$  \kb{are defined by \eqref{eq:valloss}, \eqref{eq:jacloss} and \eqref{eq:monloss} below}. 

For the loss function 
\eqref{eq:lossFunctional}, we have error in sampled values given by
\begin{equation}\label{eq:valloss}
    \mathcal{L}_{i,0}^l(\theta, U) = \dfrac{1}{n_s} \sum_{s} \left( \widetilde{\rm DtN}_{H,i}^{l}\left(\theta; u^s  \right) -  {\rm DtN}_{H,i}^{l}\left( u^s  \right) \right)^2 \,,
\end{equation}
and error with respective derivatives given by
\begin{equation}\label{eq:jacloss}
    \mathcal{L}_{i,1}^l(\theta, U) = \dfrac{1}{n_s} \sum_{k=1}^{\NHi} \sum_{s} \left( \partial_{u_k} \widetilde{\rm DtN}_{H,i}^{l}\left(\theta; u^s \right) -  \partial_{u_k}  {\rm DtN}_{H,i}^{l}\left( u^s \right) \right)^2.
\end{equation}
\kb{In addition we introduce the monotonously loss terms defined either by}
\begin{equation}\label{eq:monloss}
    \mathcal{L}_{i,\text{mon}}^l(\theta) = 
   \int_{(u_{\text{min}}, u_{\text{max}})^{\NHi}} \left( \left( \partial_{u_l} \widetilde{\rm DtN}_{H,i}^{l}\left(\theta; u  \right) \right)^{-} \right)^2 + \sum_{k \neq l} \left( \left( \partial_{u_k} \widetilde{\rm DtN}_{H,i}^{l}\left(\theta;u\right) \right)^{+} \right)^2 du, 
\end{equation}
or by
\begin{equation}\label{eq:monloss2d}
        \mathcal{L}_{i,\text{mon}}(\theta) = 
   \int_{(u_{\text{min}}, u_{\text{max}})^{\NHi}}  \left( \left( \partial_{u_l} \widetilde{\rm DtN}_{H,i}^{l}\left(\theta; u  \right) \right)^{-} \right)^2  du,
\end{equation}
where $(x)^+ = \max(x, 0)$, $(x)^- = \min(x, 0)$ for any real $x$, and \kb{where the integrals are taken over the training domain, which is an $\NHi$ dimensional cube.} 
\kb{We note that $\mathcal{L}_{i,\text{mon}}^l(\theta)$ does not require any data sampling. The purpose of this terms in the loss function is to enforce certain monotonicity properties that are expected from the original $\dtn_{H,i}$ maps. In particular, for one-dimensional problem,
we expect the Jacobian of $\dtn_{H,i}$ to have positive diagonal elements and non-positive off-diagonal ones. We note that this property is closely related to the discrete maximum principle and is classical for Finite Volume methods. The loss term defined in \eqref{eq:monloss} aims to penalize the wrong Jacobian signs. For problems in 2D, the off-diagonal elements of $\dtn'_{H,i}$ do not have specific sign in general. Therefore, for such problems, we will be using \eqref{eq:monloss2d} instead of \eqref{eq:monloss}. 
The monotonicity constrained based on \eqref{eq:monloss} or \eqref{eq:monloss} allows us to improve the robustness of Newton's method for \eqref{eq:approxsubstlearn}, which otherwise would be likely to fail.
}

For training, we  generate 
${\rm DtN_{H,i}}(u^s)$ and  ${\rm DtN'_{H,i}}(u^s)$ values for each sampling vector $(u^s)_{s=1}^{n_s}$.
%
%
%
%
In terms of choosing sampling points to be used in training, 
 we use a uniform sampling procedure.%
 Given a sampling range $(u_{min}, u_{max})$, we sample each vector component with $m$ entries equally spaced between $u_{min}$ and $ u_{max}$. Occasionally, these $m$ sampling points can also include additional points, particularly in the ``center" of existing training points. This naturally leads to $n_s= m^{\NHi}$ training vectors for the NN model,
 where $(u^s) \in \mathbb{R}^{\NHi}$.

\section{Numerical experiments}\label{sec:numericalExp}

\subsection{1D degenerate elliptic problem}
\begin{figure}[h]
\centering
	\centering
	\includegraphics[width=.343\linewidth]{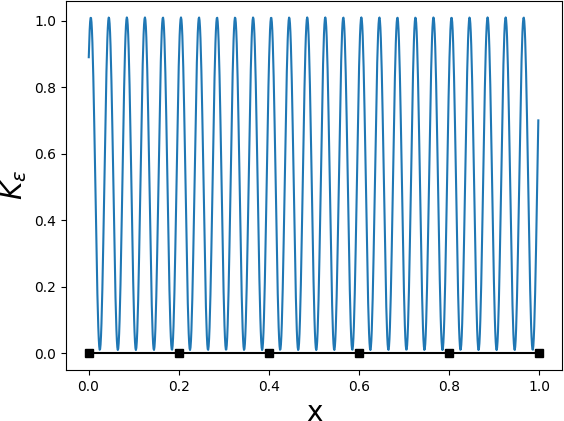}
\caption{Partitioning of the model domain and the plot of $K_\epsilon(x)$}\label{fig_pme_K}
\end{figure}

We first consider the following nonlinear equation 
\begin{equation}\label{pme}
\mass u - \left( K_\epsilon( x ) \left( u^4 \right)' \right)' = 0 
\end{equation}
posed in the domain $(0, 1)$ with some non-negative Dirichlet boundary conditions. The diffusion coefficient is given by 
$$
K_\epsilon( x ) = 10^{-2} + 
\frac{1}{2}\left( 1 + \sin\left( 10 \pi x + \frac{\pi}{4} \right)  \right),
$$
and $\mass=20$. Equation \eqref{pme} can be interpreted as a stationary variant of the porous media equation \cite{vazquez2007porous}. It can be shown that $u$ takes values in $[0, \max(u(0), u(1)]$. \kb{We note that, since} the diffusion term in \eqref{pme} degenerates at $u=0$, the solutions of \eqref{pme} may vanish over some portion of the domain. The domain $(0,1)$ is partitioned in five subdomains of equal size. We refer to Figure \ref{fig_pme_K} for the illustration of the partitioning and the coefficient $K_\epsilon$. 
Because the partitioning of the domain is matching periodicity of $K_\epsilon(x)$, the Dirichlet-to-Neumann operators coincide for all of the subdomains.  As 
the analytical expression of the subdomain's Dirichlet-to-Neumann operator is not available, the latter will be replaced by the approximate one \kb{denoted by $\dtn_{h,i}$} and computed by a finite element method using $200$ grid points. 
To be more specific, we use $\mathbb{P}_1$ finite element method with mass lumping. \kb{We note that in the one-dimensional case the mapping $\dtn_{h,i}$, acting from $\mathbb{R}^2$ to itself, coincides with the coarse operator $\dtn_{H,i}$}
For surrogate model $\wdtn_{H,i}$ we use fully connected neural network with two hidden layers of 64 neurons each. \kb{The activation function used in all the numerical experiment is given by $\sigma(u) = \max(u,0)^2$.} The model is trained over the domain $[0, u_{max}]^2$ with $u_{max} = 4$. Training points are sampled over a regular grid in $[0, u_{max}]^2$ consisting of $n_s = 2^2, 3^2, 4^2$ and $5^2$ points.

The loss function involves both values and derivatives of $\dtn_{H,i}$ at the sampling points with $c_0 = 1$ and $c_1 = 0.1$, as well as the monotonicity component \eqref{eq:monloss} with $c_{mon} = 4$. The integral in \eqref{eq:monloss} is approximated numerically based on integrand values over a regular $40\times 40$ grid.

\begin{figure}[h]
\centering
\begin{subfigure}{0.33\textwidth}
	\centering
	\includegraphics[height=5cm, width=.9\linewidth]{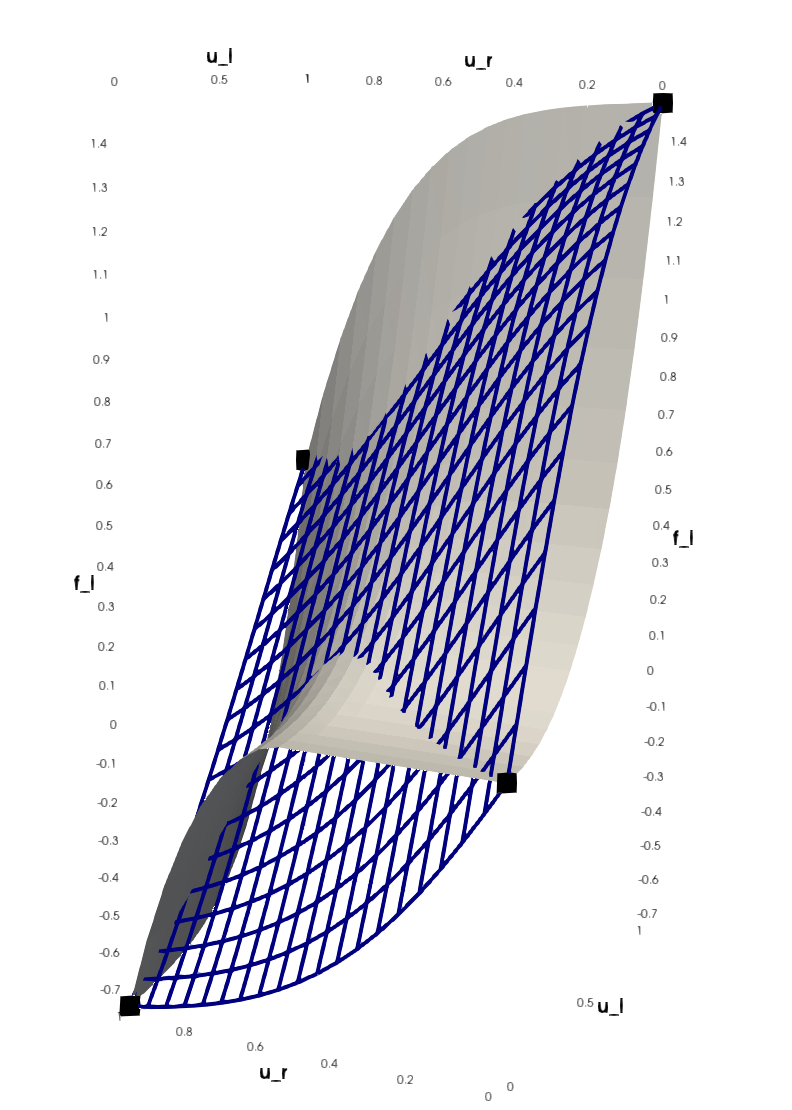}
\end{subfigure}
\begin{subfigure}{0.33\textwidth}
	\centering
	\includegraphics[height=5cm, width=.85\linewidth]{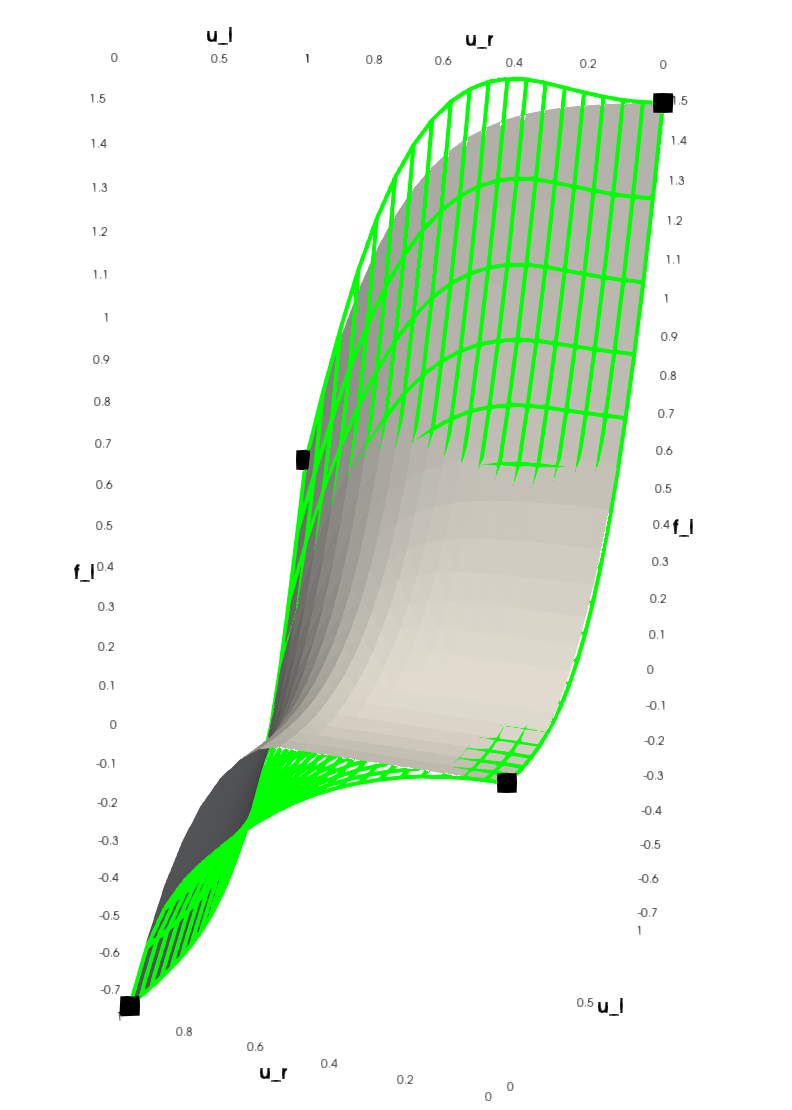}
\end{subfigure}
\begin{subfigure}{0.33\textwidth}
	\centering
	\includegraphics[height=5cm, width=.85\linewidth]{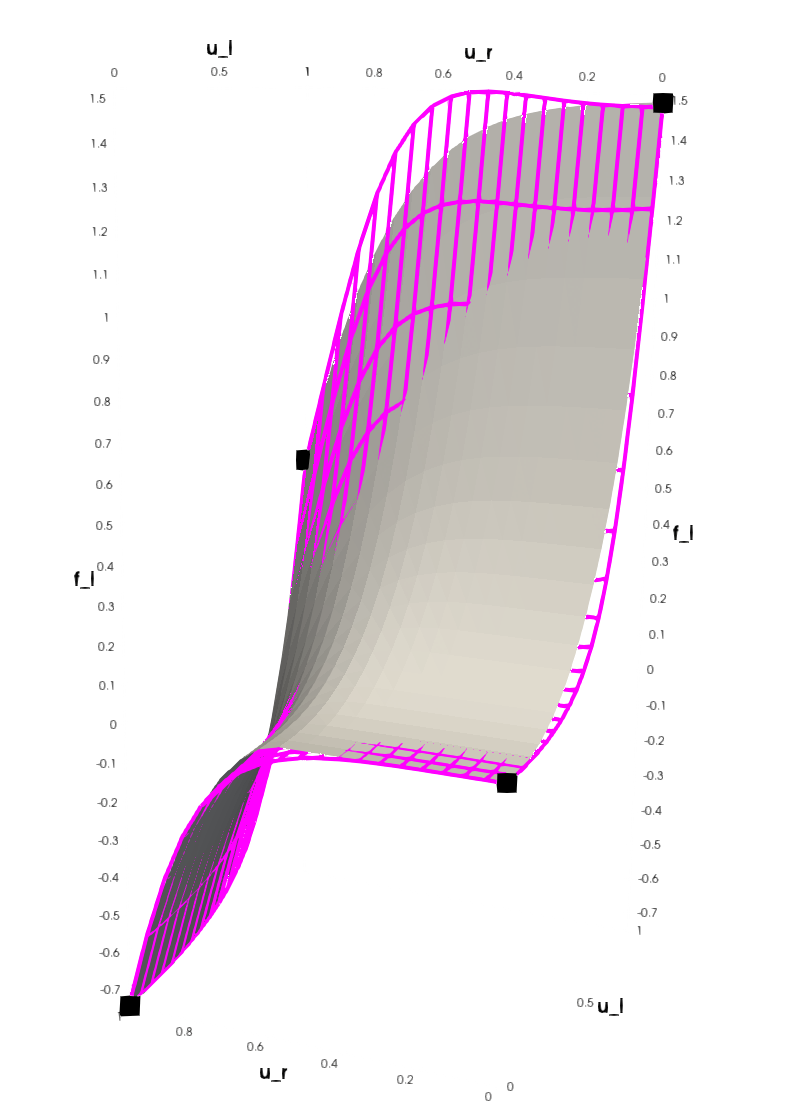}
\end{subfigure}
\newline 
\begin{subfigure}{0.33\textwidth}
	\centering
	\includegraphics[height=5cm, width=.9\linewidth]{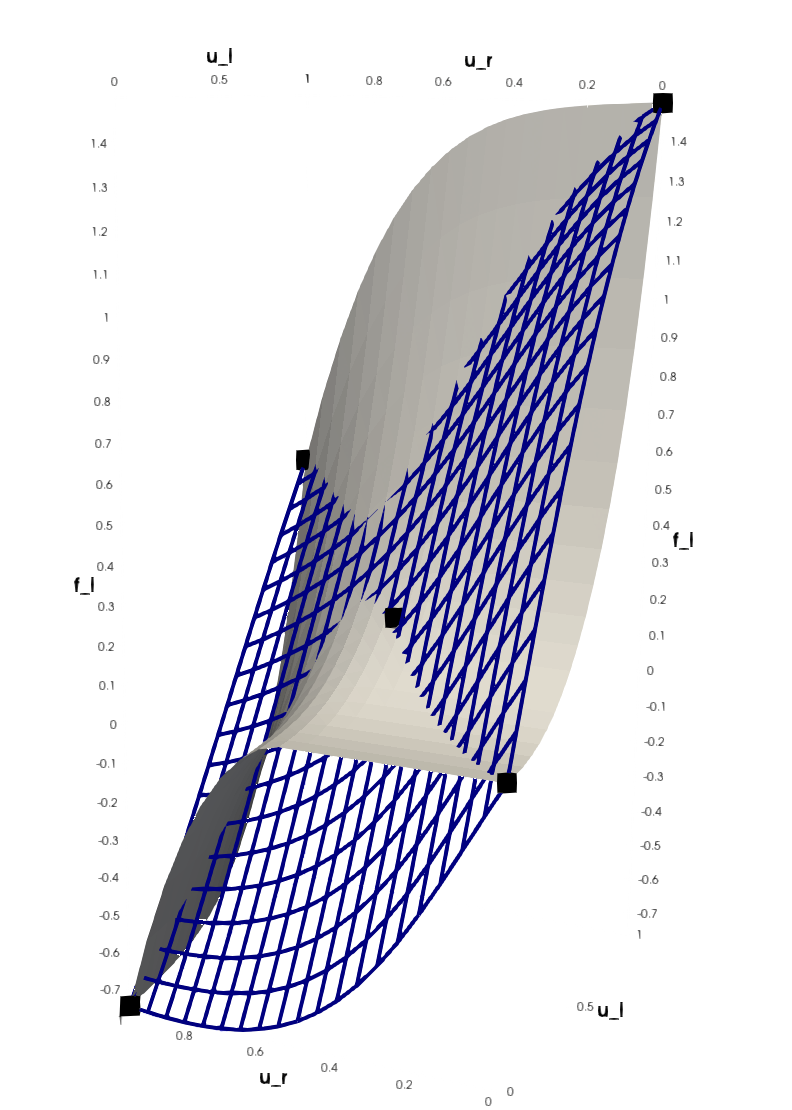}
\end{subfigure}
\begin{subfigure}{0.33\textwidth}
	\centering
	\includegraphics[height=5cm, width=.85\linewidth]{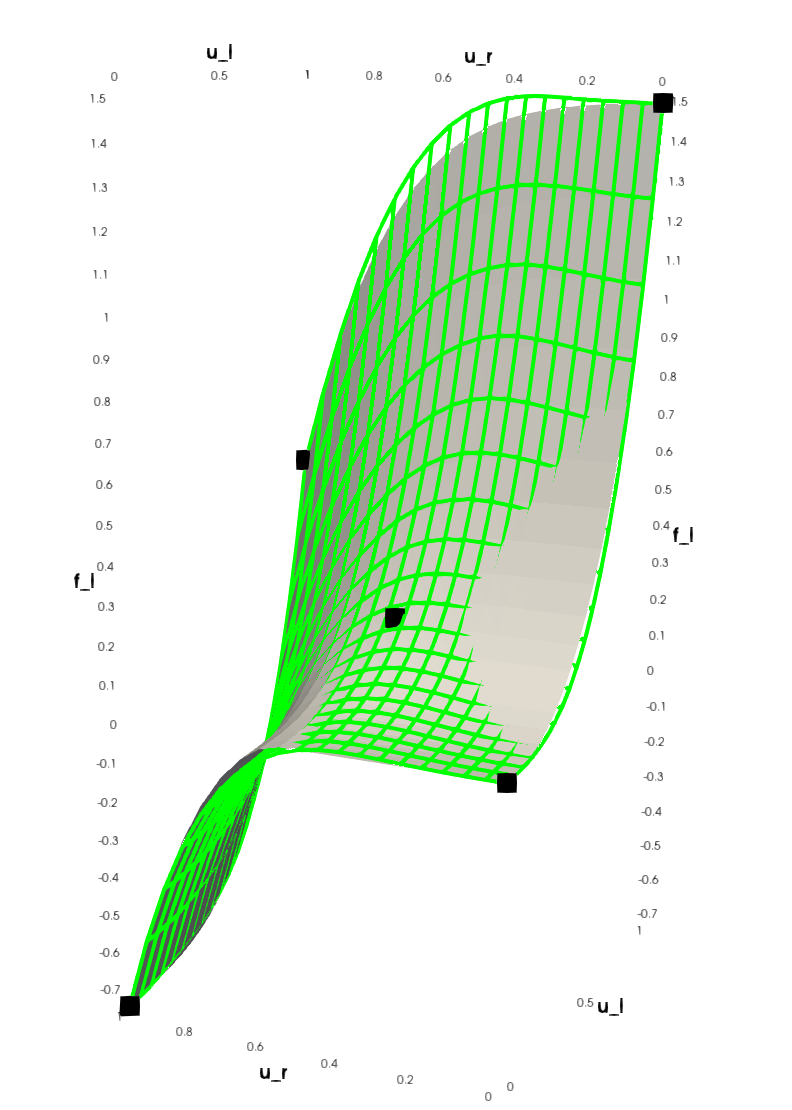}
\end{subfigure}
\begin{subfigure}{0.33\textwidth}
	\centering
	\includegraphics[height=5cm, width=.85\linewidth]{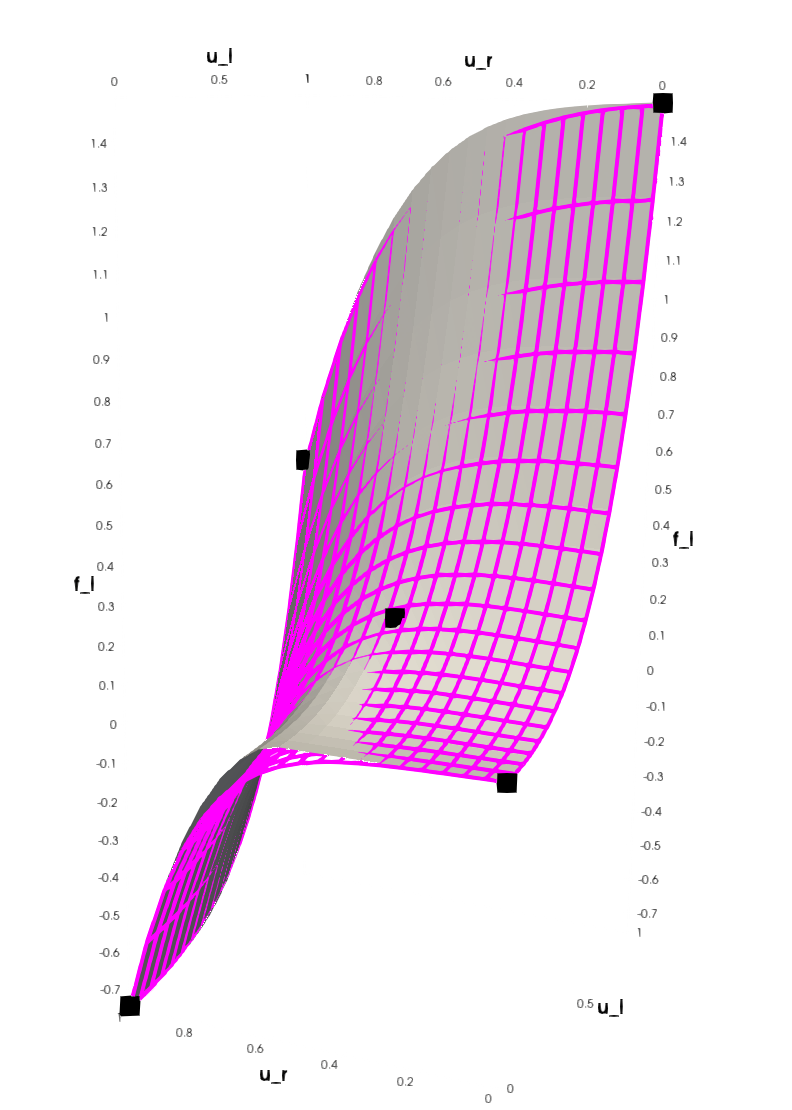}
\end{subfigure}
\newline 
\begin{subfigure}{0.33\textwidth}
	\centering
	\includegraphics[height=5cm, width=.85\linewidth]{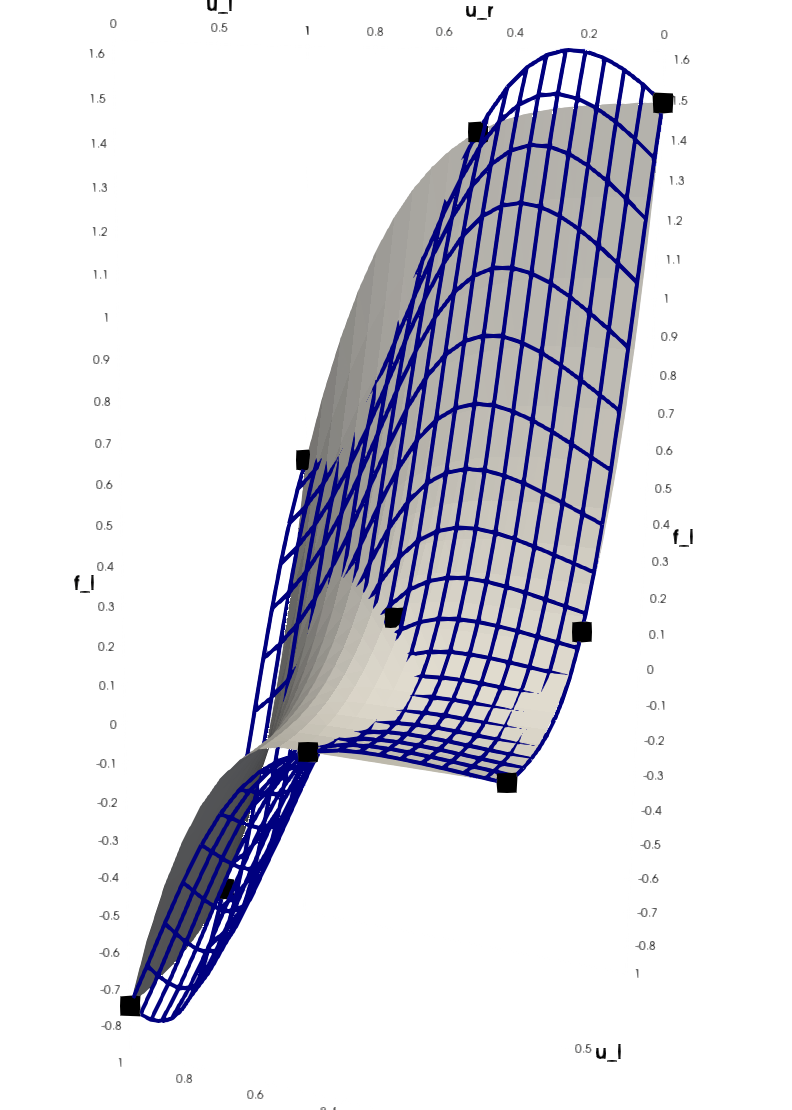}
\end{subfigure}
\begin{subfigure}{0.33\textwidth}
	\centering
	\includegraphics[height=5cm, width=.85\linewidth]{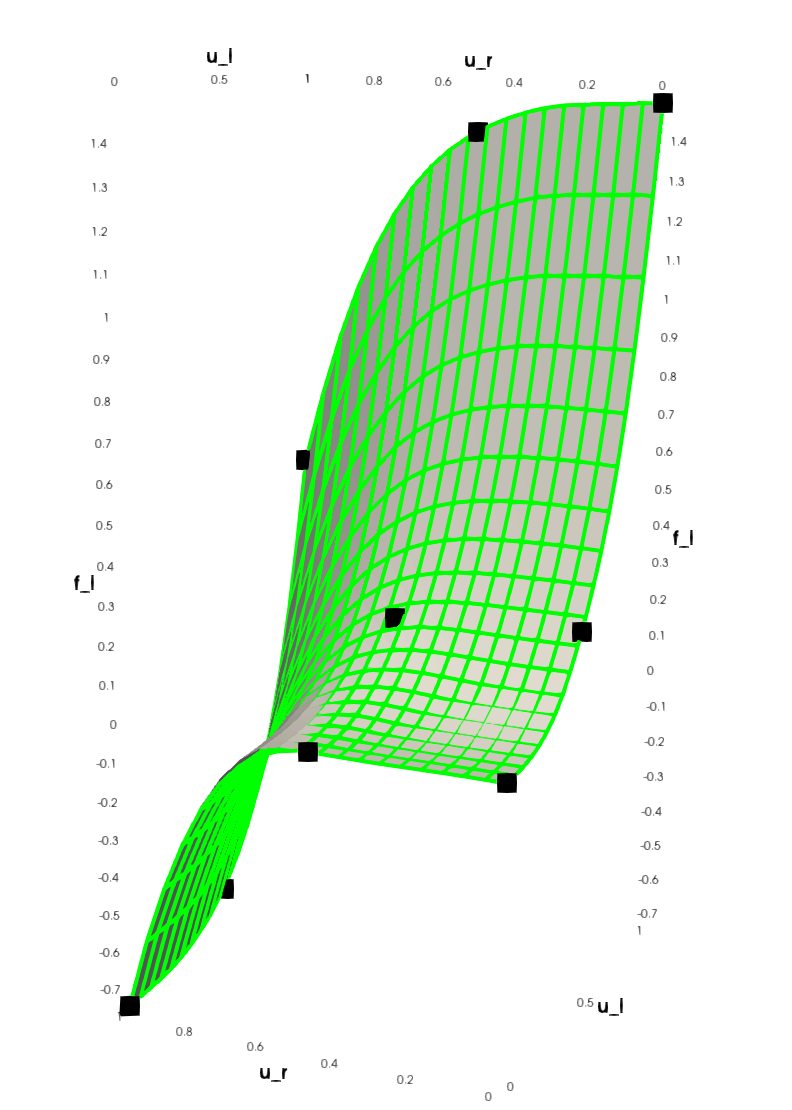}
\end{subfigure}
\begin{subfigure}{0.33\textwidth}
	\centering
	\includegraphics[height=5cm, width=.85\linewidth]{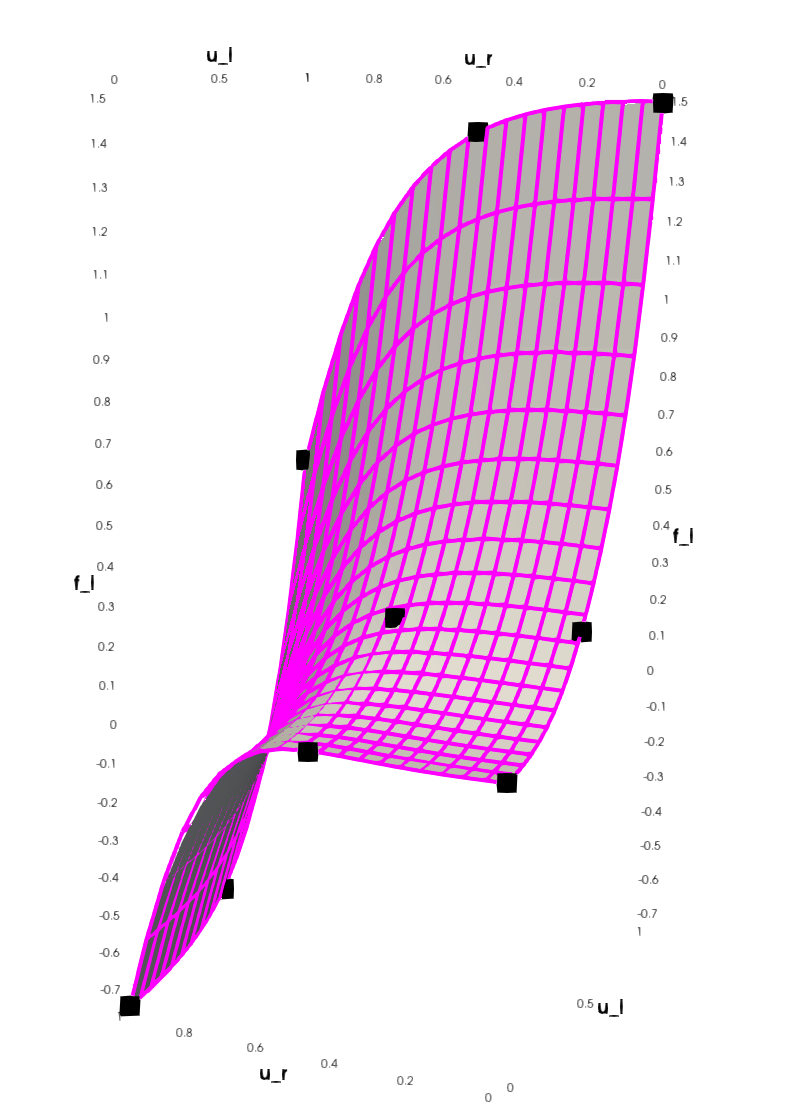}
\end{subfigure}
\caption{
Qualitative comparison of the surfaces $z = \dtn_{H,i}((u_{\rm left}, u_{\rm right}))^{1}$ (solid) and $z = \wdtn_{H,i}((u_{\rm left}, u_{\rm right}))^{1}$ (wireframe). Left column (blue): only DtN values included. Middle column (green): DtN and DtN derivatives included. Right column (purple): DtN values, derivatives, and monotonicity assertion included. Top to bottom: $n_s = 2^2$ training points, $n_s=2^2+1$ training points, $n_s=3^2$ training points.}
\label{fig:learneddtn}
\end{figure}

\textcolor{blue}{Let us begin with a qualitative analysis of the surrogate $\wdtn_{H,i}$ model. Recall that, in 1D, $\dtn_{H,i}$ (coinciding with $\dtn_{h,i}$) maps $\mathbb{R}^2$ to itself, which can be interpreted as mapping left and right Dirichlet data to a pair of left and right fluxes. We report on Figure \ref{fig:learneddtn}  the left flux and its surrogates as the function of the input Dirichlet values. The surface
$z = \dtn_{H,i}((u_{\rm left}, u_{\rm right}))^{1}$ shown as solid is the same for all sub-figures. Colored wireframe surfaces show $z = \wdtn_{H,i}((u_{\rm left}, u_{\rm right}))^{1}$ for various number of sampling points (4, 5 and  9 for the first, second and third rows respectively), and different values of coefficients $c_\alpha$ in \eqref{eq:lossFunctional} (same for each column). The results obtained using $\dtn$ values alone ($c_1 = 0 = c_{mon} = 0$) are shown in the first column, those obtained using both values and derivatives ($c_{mon} = 0$) are reported in the middle column. Finally, the right column shows the outputs of the model using the additional monotonicity loss term \eqref{eq:monloss}.
For the sake of clearer visualisation, here,the  parameter $a$ is set to 1.
Unsurprisingly, fitting not only the values, but also the tangential plains has a drastic benefit on the quality of the surrogate model (compare left and middle columns). Models with $c_1 > 0$ manage to capture the shape of the response surface reasonably well for very few training points. We also observe that adding the monotonicity loss does not seems to pollute the interpolation quality (middle versus right columns).}

\begin{figure}[h]
\centering
\begin{subfigure}{0.49\textwidth}
	\centering
	\includegraphics[width=.7\linewidth]{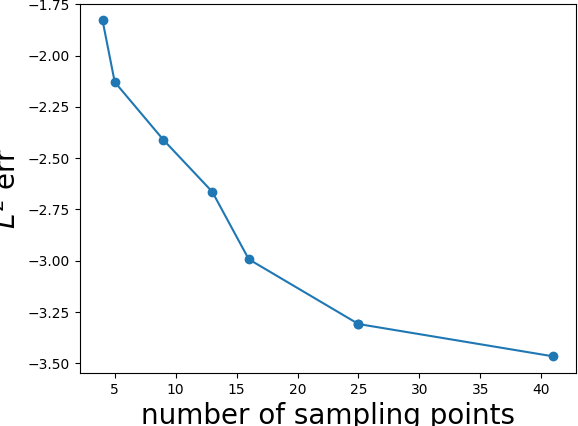}
\end{subfigure}
\begin{subfigure}{0.49\textwidth}
	\centering
	\includegraphics[width=.7\linewidth]{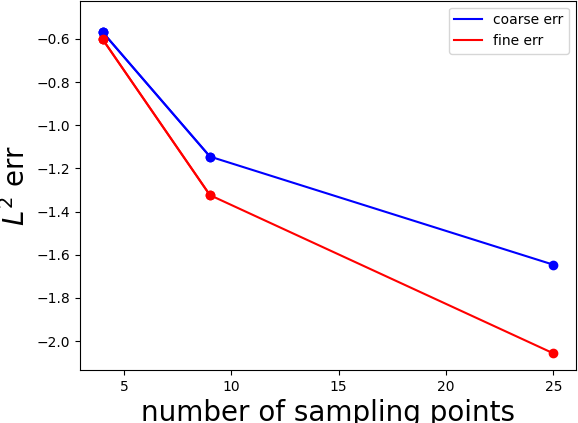}
\end{subfigure}
\caption{Left: Interpolation relative $L^2$ error. Right: Solution relative $L^2$ error.}\label{fig_pme_conv}
\end{figure}

We report on Figure \ref{fig_pme_conv} convergence of the interpolation error of $\dtn_{H,i}$ in relative $L^2( (0, u_{max})^2 )$ norm as the function of the number of sampling points. The error is evaluated over a regular $20 \times 20 $ grid in $[0, u_{max}]^2$. \kb{Here, in addition to the points sampled over regular $n_s\times n_s$ grid, we include additional points  taken at the centers of the sampling grid cells.}
The decrease of error deteriorates after approximately 25 sampling points. This saturation effect (around error value of $10^{-3.5}$) is even more clear in the second test case (see Figure \ref{fig_plap_conv} below). The failure in achieving higher accuracy should probably be attributed to our current training algorithm, combined with the fact that the structure of the surrogate model is fixed through the whole experiment.

\begin{figure}[h]
\centering
\begin{subfigure}{0.49\textwidth}
	\centering
	\includegraphics[width=.7\linewidth]{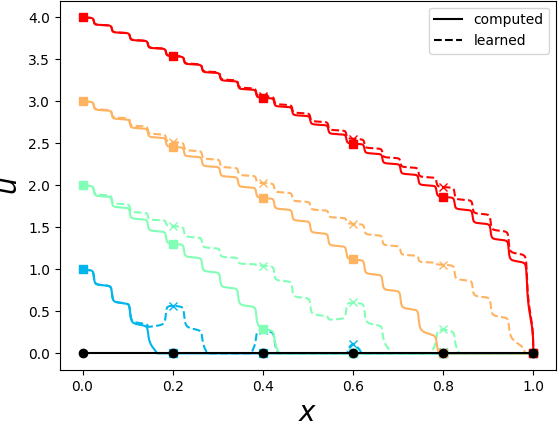}
\end{subfigure}
\begin{subfigure}{0.49\textwidth}
	\centering
	\includegraphics[width=.7\linewidth]{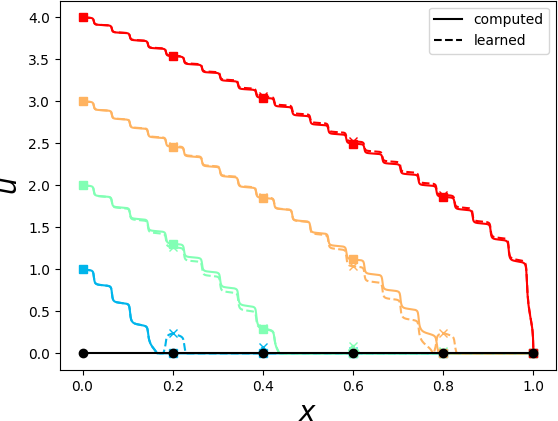}
\end{subfigure}
\begin{subfigure}{0.49\textwidth}
	\centering
	\includegraphics[width=.7\linewidth]{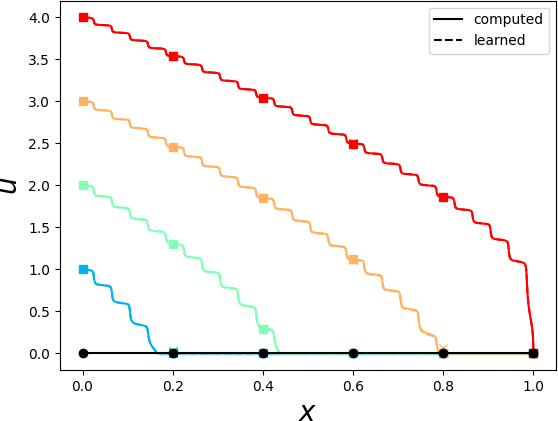}
\end{subfigure}
\begin{subfigure}{0.49\textwidth}
	\centering
	\includegraphics[width=.7\linewidth]{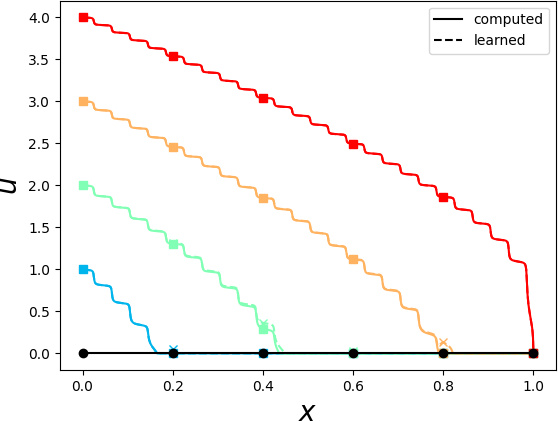}
\end{subfigure}
\caption{Solution profiles \mb{of \eqref{pme} for various left boundary conditions} using $\dtn_{H,i}$ (solid) and $\wdtn_{H,i}$ (dashed) operators for $n_s = 4, 9, 16$ and $25$.} \label{fig_pme_surrogate}
\end{figure}

We further evaluate the performance of the surrogate models for the solution of the substructured problem \eqref{eq:approxsubstlearn}, with $u(1) = 0$, and $u(0)$ taking values in $\{1, 2, 3, 4\}$. \kb{Once the coarse solution $g_H$ and $\widetilde{g}_H$ are computed by approximately 
solving \eqref{eq:nonlincoarsesubstruct}, and \eqref{eq:approxsubstlearn} respectively, we reconstruct the solution inside the subdomains by solving local problems on $\Oi$ with boundary condition given either by $g_H|_{\partial\Oi}$ or $\widetilde{g}_H|_{\partial\Oi}$. }
Figure \ref{fig_pme_surrogate} reports \kb{the reconstructed solution, using both the original and the surrogate $\dtn$ models,} for different boundary conditions and varying number of sampling points. Starting from $n_s = 3^2$ the surrogate model produced the solutions which are visually appear to be accurate.  
The relative $L^2(0,1)$ error as the function of the number of sampling points is reported by figure \ref{fig_pme_conv}. 

\begin{figure}[h]
	\centering
	\includegraphics[width=.343\linewidth]{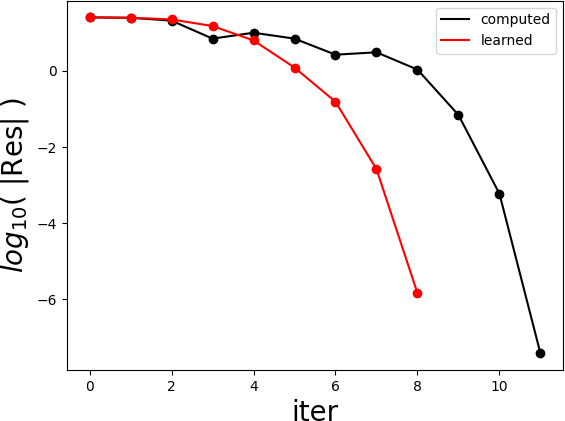}
\caption{Convergence of Newton's method using $\dtn_{H,i}$ (black) and $\wdtn_{H,i}$ (red) operators to solve \eqref{pme}.} \label{fig_pme_newton}
\end{figure}

The numerical solution of \kb{\eqref{eq:nonlincoarsesubstruct} and \eqref{eq:approxsubstlearn}} is obtained by Newton's method. On Figure \ref{fig_pme_newton}, we report typical convergence history for Newton's method for the original \eqref{eq:nonlincoarsesubstruct} and the surrogate model \eqref{eq:approxsubstlearn} with $u(0) = u_{max}$ and $n_s = 4^2$. Here, the use of the surrogate model leads to convergence of Newton's method with slightly fewer iterations. Most importantly, the inference time of $\wdtn$ is much lower than the time required to evaluate $\dtn$. 

\subsection{1D $p-$Laplace problem}
\begin{figure}[h]
\centering
\begin{subfigure}{0.49\textwidth}
	\centering
	\includegraphics[width=.7\linewidth]{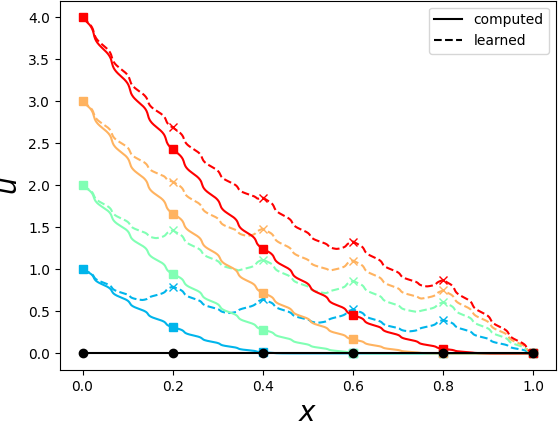}
\end{subfigure}
\begin{subfigure}{0.49\textwidth}
	\centering
	\includegraphics[width=.7\linewidth]{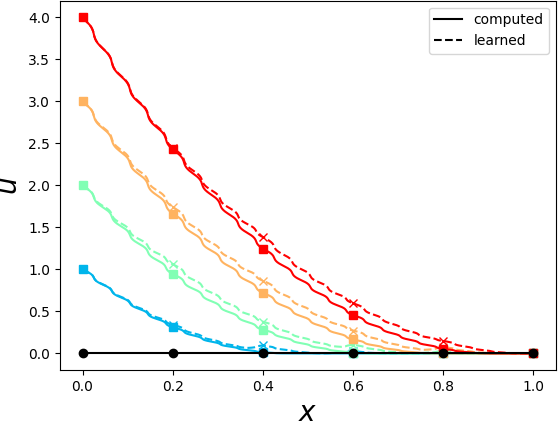}
\end{subfigure}
\begin{subfigure}{0.49\textwidth}
	\centering
	\includegraphics[width=.7\linewidth]{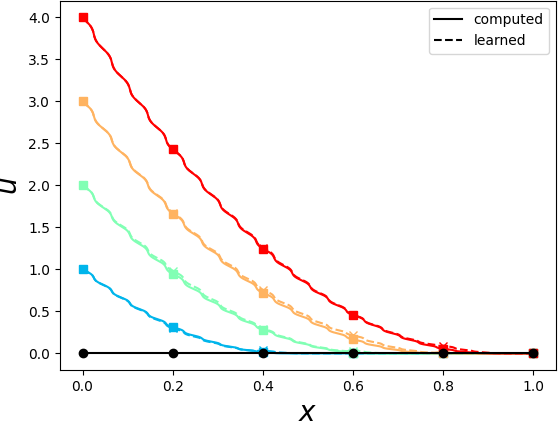}
\end{subfigure}
\begin{subfigure}{0.49\textwidth}
	\centering
	\includegraphics[width=.7\linewidth]{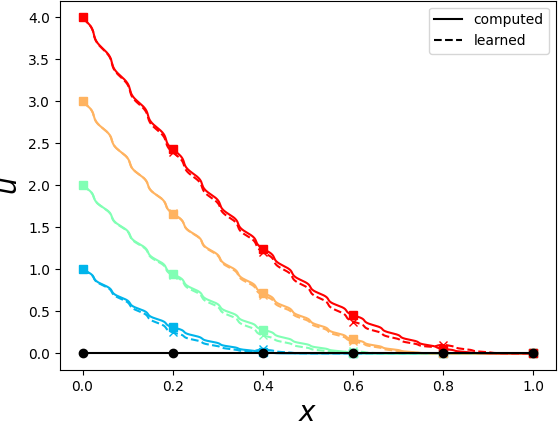}
\end{subfigure}
\caption{Solution profiles  \mb{of \eqref{pme} for various left boundary conditions} using $\dtn_{H,i}$ (solid) and $\wdtn_{H,i}$ (dashed) operators for $n_s = 4, 9, 16$ and $25$} \label{fig_plap_surrogate}
\end{figure}

\begin{figure}[h]
\centering
\begin{subfigure}{0.49\textwidth}
	\centering
	\includegraphics[width=.7\linewidth]{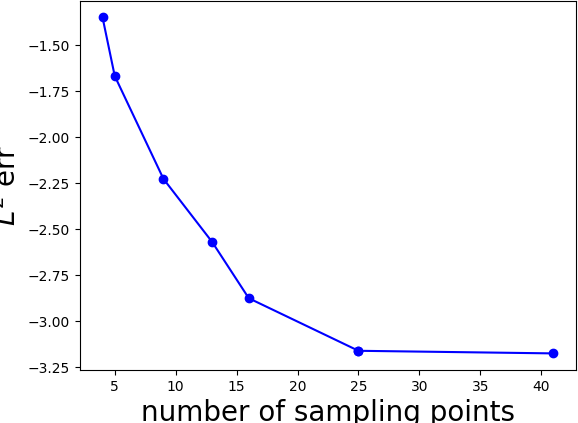}
\end{subfigure}
\begin{subfigure}{0.49\textwidth}
	\centering
	\includegraphics[width=.7\linewidth]{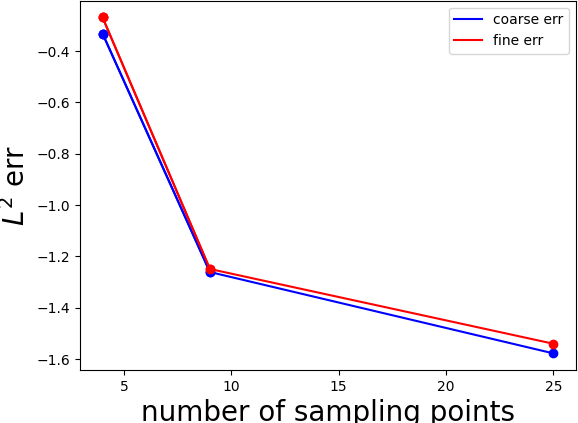}
\end{subfigure}
\caption{Left: Interpolation relative $L^2$ error. Right: Solution relative $L^2$ error.}\label{fig_plap_conv}
\end{figure}

\begin{figure}[h]
	\centering
	\includegraphics[width=.343\linewidth]{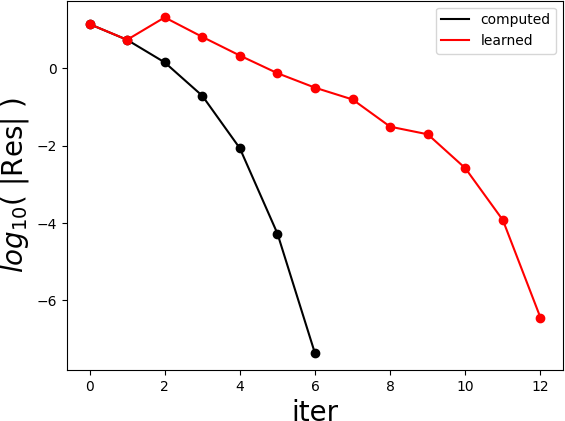}
\caption{Convergence of Newton's method using $\dtn$ (black) and $\wdtn$ (red) operators to solve \eqref{plap}.}
\label{fig_plap_newton}
\end{figure}

Next, we consider the following equation 
\begin{equation}\label{plap}
\mass u - \left( K_\epsilon( x ) | u' |^2 u' \right)' = 0
\end{equation}
\kb{All the numerical parameters are set up as in the previous test case, except for the coefficient $\mass = 5$. Again we compare the solution obtained using either the original or the surrogate substructured problem (\eqref{eq:nonlincoarsesubstruct} or \eqref{eq:approxsubstlearn}). Solution comparison for different boundary conditions and varying sampling point size is reported reported on Figure \ref{fig_plap_surrogate}. Starting from $n_s = 3^2$ the solution obtained using the surrogate model appears to be reasonably accurate.  Figure \ref{fig_plap_conv} exhibits convergence of the relative interpolation and solution errors in $L^2((0,u_{max})^2)$ and $L^2(0,1)$ norms,respectively. As before, we observe the saturation of the interpolation error as the function of $n_s$; this phenomenon  is even more pronounced than in the previous test case.} Typical convergence of Newton's method is shown by Figure \ref{fig_plap_newton} for $u(0) = 4$ and $n_s = 4^2$. Although Newton's method requires twice as much iterations for the surrogate model compared to the original one, the former leads to a much faster algorithm due to a much lower inference time.

\subsection{2D degenerate elliptic problem}
\begin{figure}[h]
\centering
\begin{subfigure}{0.49\textwidth}
	\centering
	\includegraphics[width=.7\linewidth]{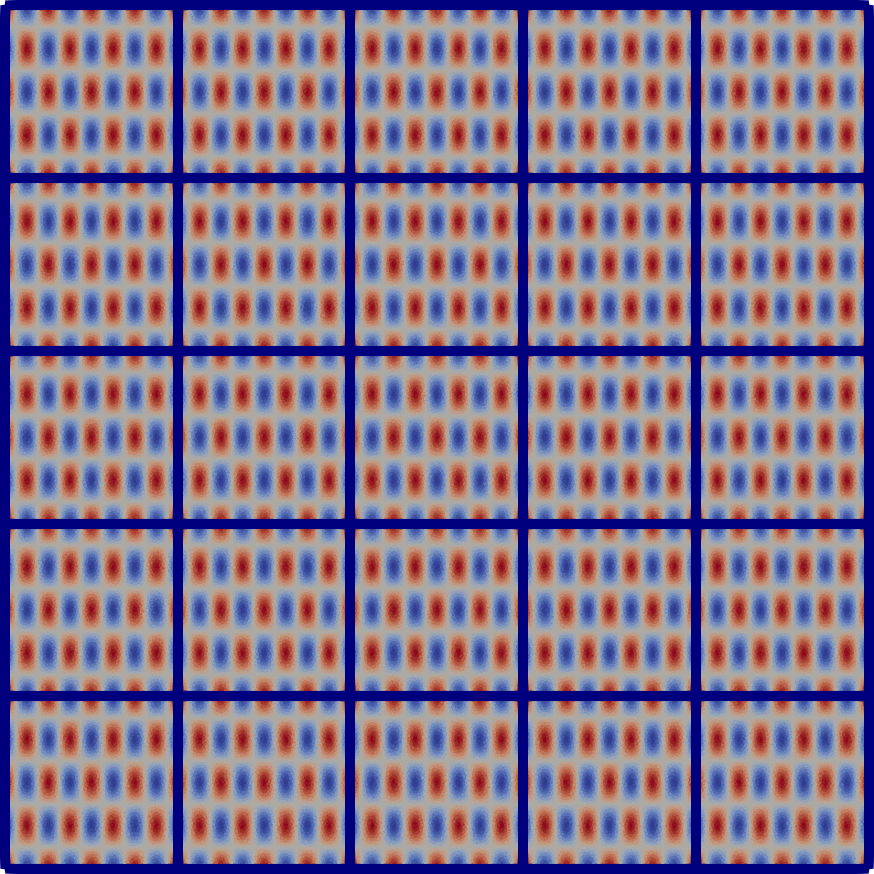}
\end{subfigure}
\begin{subfigure}{0.49\textwidth}
	\centering
	\includegraphics[width=.7\linewidth]{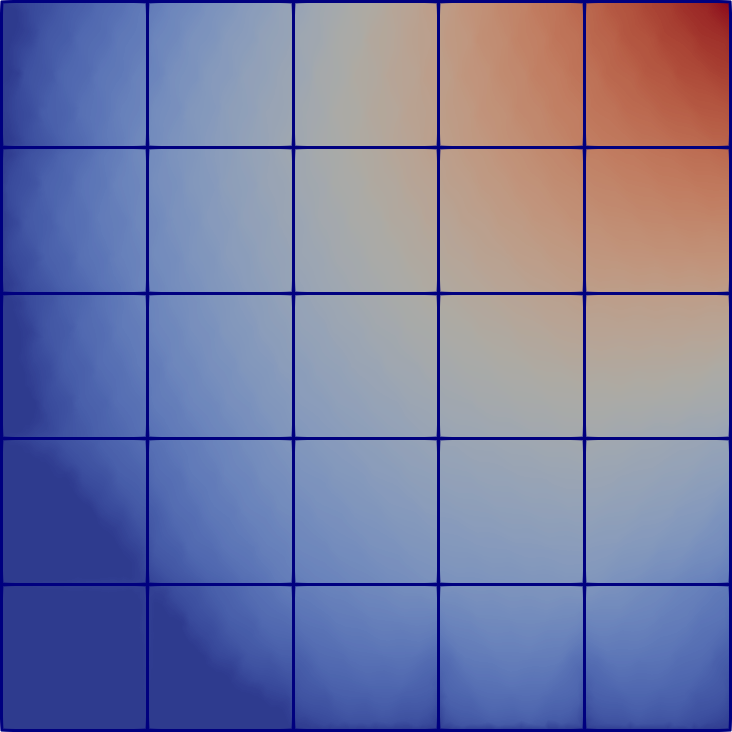}
\end{subfigure}
\caption{Left: Diffusion coefficient $K_\epsilon(x)$. Right: Solution of model problem \eqref{pme2d}  using the original $\dtn_{H,i}$ operator. 
}\label{fig_pme_K_2d}
\end{figure}

We consider a bi-dimensional version of \eqref{pme} given by 
\begin{equation}\label{pme2d}
\begin{array}{rlllll}
\dsp u - {\rm div}\left( K_\epsilon(x,y) \nabla u^4 \right) =0 & \mbox{in} & \Omega = (0,1)^2, \\
\end{array}
\end{equation}
combined with the Dirichlet boundary conditions 
$$
u(x,y) = \max( u_{max}( x + y - 1 ), 0 ) \qquad \mbox{on}  \qquad \partial \Omega
$$
where $u_{max} = 1.2$. That is, $u(x,y) = 0$ on lower and left boundaries and $u(x,y) = u_{max}( x + y - 1 )$ on right and top ones. We set 
$$
K_{\epsilon}(x,y) =  10^{-2} + \frac{1}{2}\left( 1 + \sin\left( 10 x + \frac{\pi}{2} \right) \sin\left( 5 y + \frac{\pi}{2} \right) \right),
$$
and we consider a regular partitioning of $\Omega$ into $5\times 5$ subdomains (see Figure \ref{fig_pme_K_2d} for illustration). Here again, since the partitioning of the domain matches the periodicity of the coefficient $K_\epsilon$, one only needs to build a single surrogate $\dtn$ model. More precisely, for every component of $\dtn_{H,i}^l$, we train a specific model $\wdtn_{H,i}^l$ using a fully connected neural network with 2 hidden layers of $20$ neurons each \kb{and the activation function $\sigma(u) = \max(u,0)^2$}.  


The reference $\dtn_{H,i}$ operator is computed using mass lumped $\mathbb{P}_1$ finite element method on a fine grid with approximately 660 triangles and 340 nodal degrees of freedom by subdomain. We report on Figure \ref{fig_pme_K_2d} the solution using this reference map. 
The training data is sampled on a regular grid over $[0, u_{max}]^4$. The loss function uses both values and derivatives of $\dtn_{H,i}^l$ at the sampling points with $c_0 = 1$ and $c_1 = 0.1$. The monotonicity of $\wdtn_{H,i}^l$ is enhanced by adding the loss term \eqref{eq:monloss2d} with $c_{mon} = 10$. 
The integral in \eqref{eq:monloss2d} is approximated by Monte Carlo method, using $200$ random points sampled in $[0, u_{max}]^4$ at every optimisation step.

\begin{figure}[h]
\centering
\begin{subfigure}{0.49\textwidth}
	\centering
	\includegraphics[width=.6\linewidth]{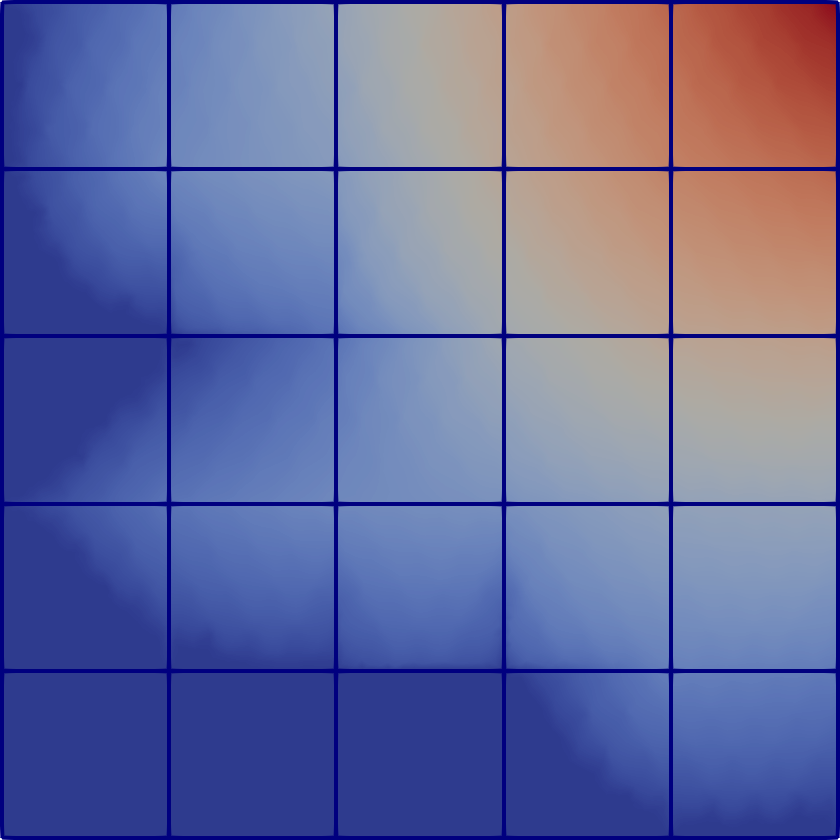}
\end{subfigure}
\begin{subfigure}{0.49\textwidth}
	\centering
	\includegraphics[width=.6\linewidth]{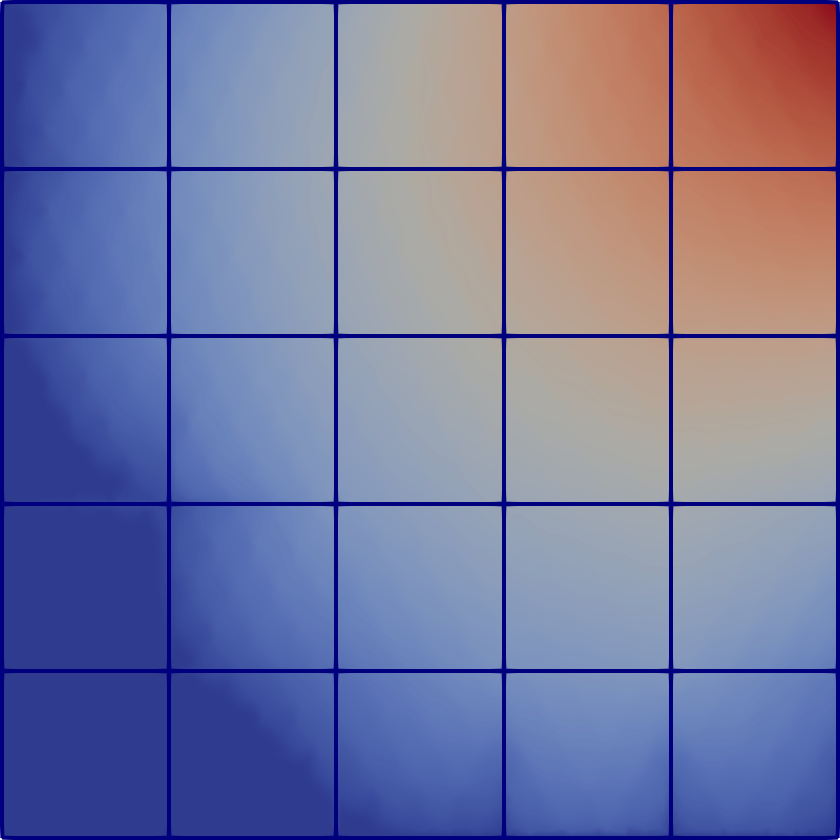}
\end{subfigure}
\begin{subfigure}{0.49\textwidth}
	\centering
	\includegraphics[width=.6\linewidth]{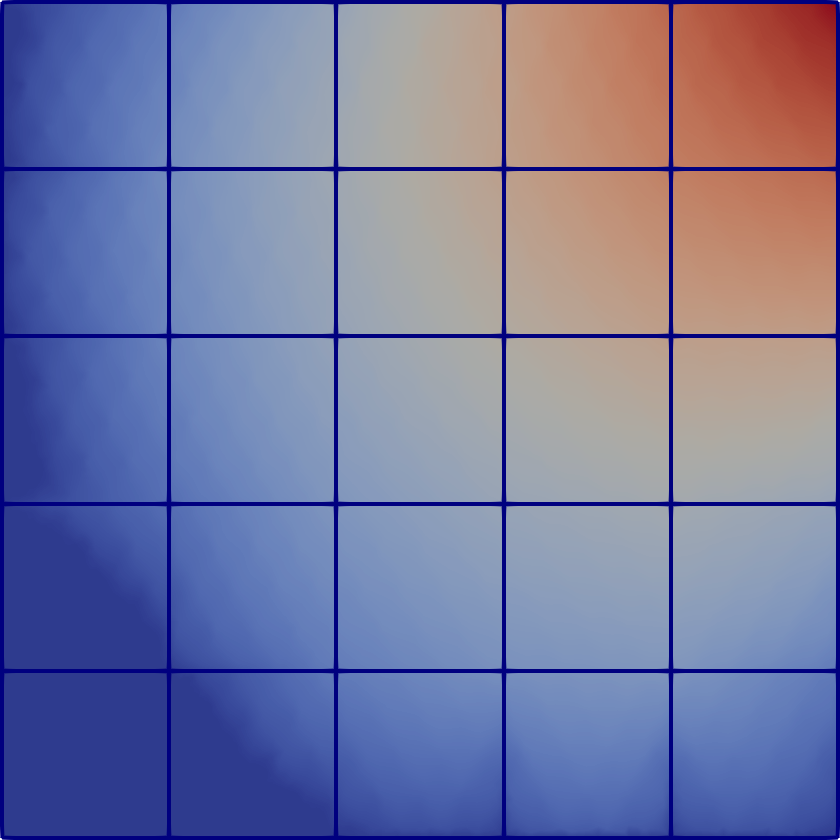}
\end{subfigure}
\begin{subfigure}{0.49\textwidth}
	\centering
	\includegraphics[width=.6\linewidth]{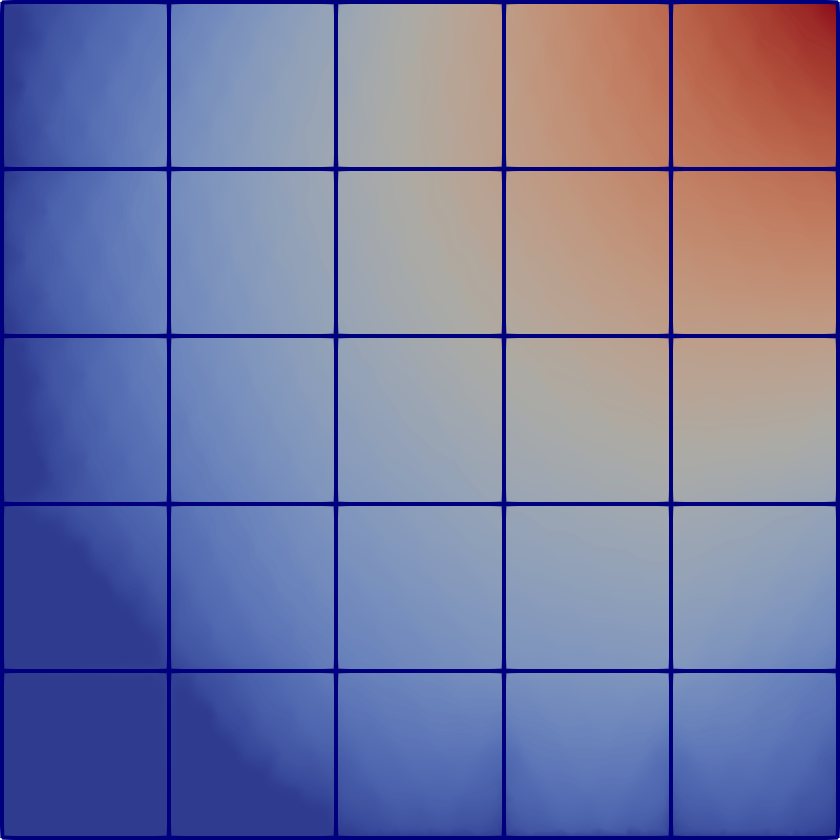}
\end{subfigure}
\caption{Solution of model problem \eqref{pme2d} using learned $\wdtn$ map for various numbers of sampling points, including $n_s = 2^4$ (top left), $n_s=3^4$ (top right), $n_s=4^4$ (bottom left),  and $n_s=5^4$ (bottom right).} \label{fig_pme_2d_surrogate}
\end{figure}

\begin{figure}[h]
\centering
\begin{subfigure}{0.49\textwidth}
	\centering
	\includegraphics[width=.7\linewidth]{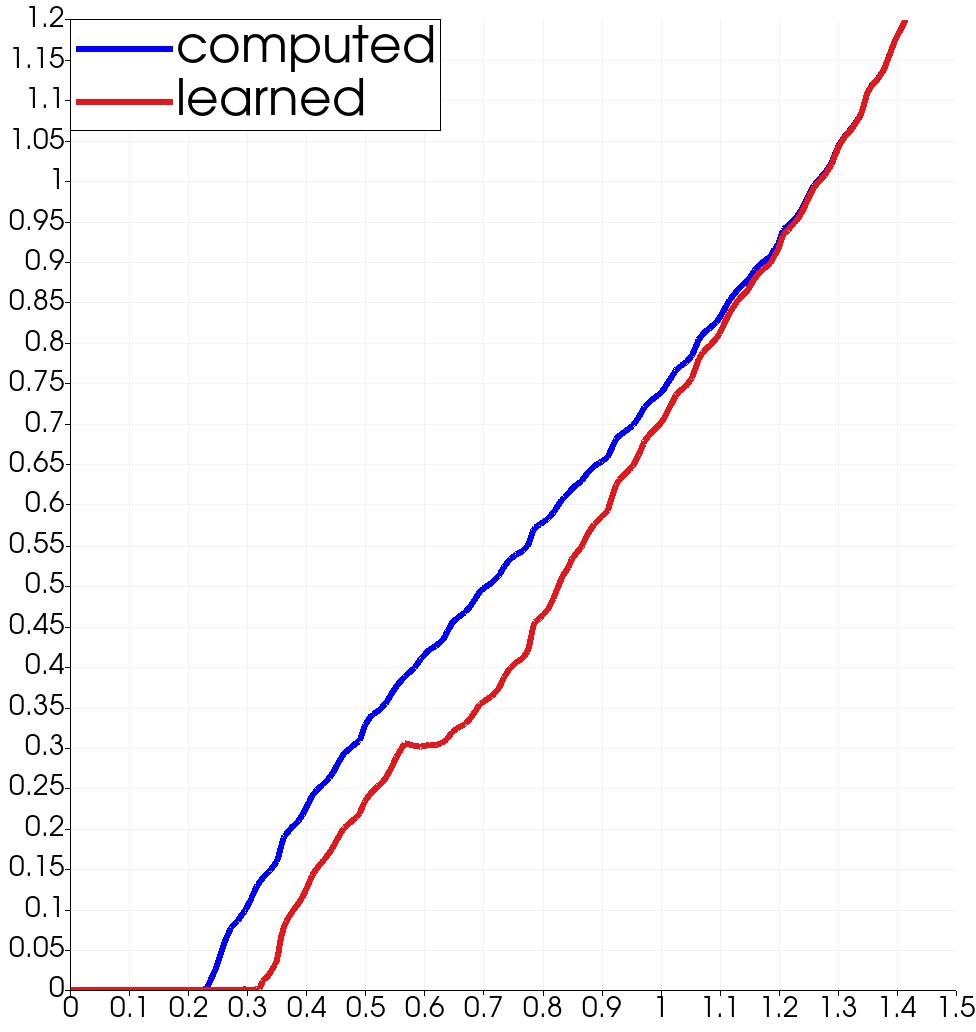}
\end{subfigure}
\begin{subfigure}{0.49\textwidth}
	\centering
	\includegraphics[width=.7\linewidth]{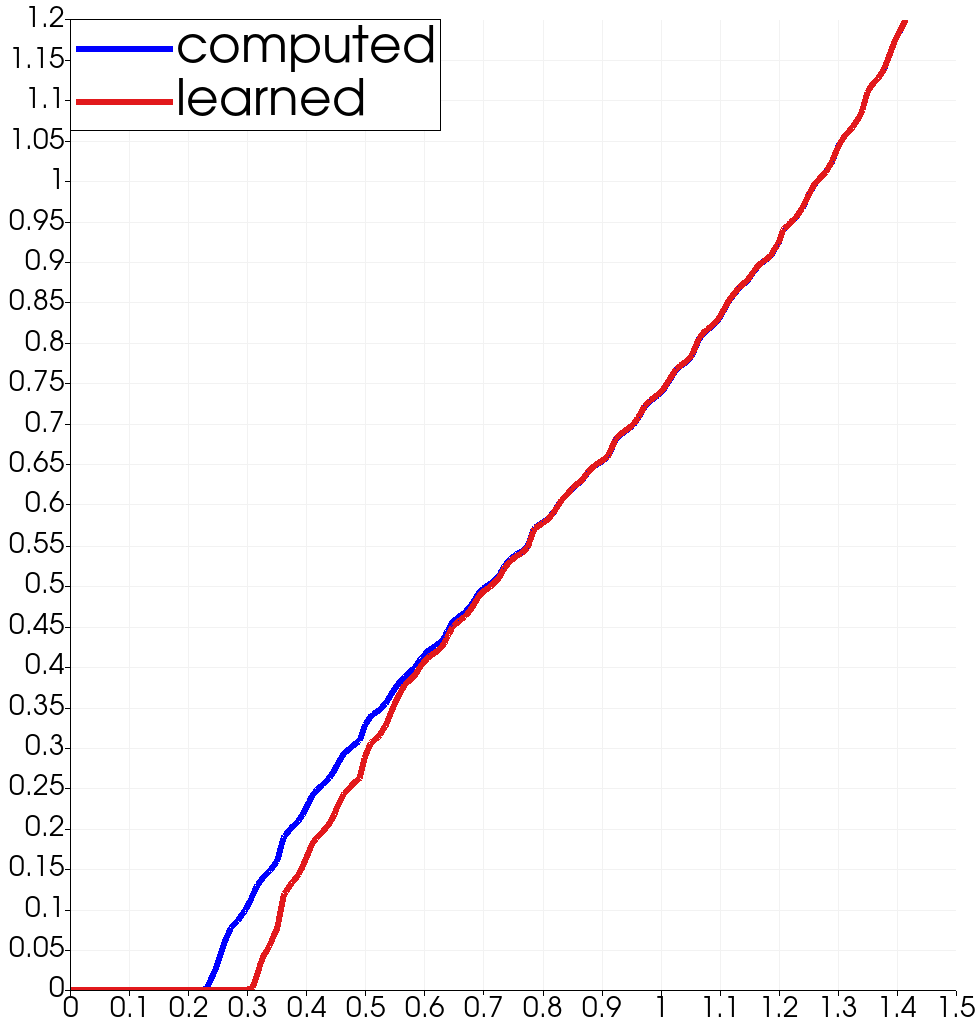}
\end{subfigure}
\begin{subfigure}{0.49\textwidth}
	\centering
	\includegraphics[width=.7\linewidth]{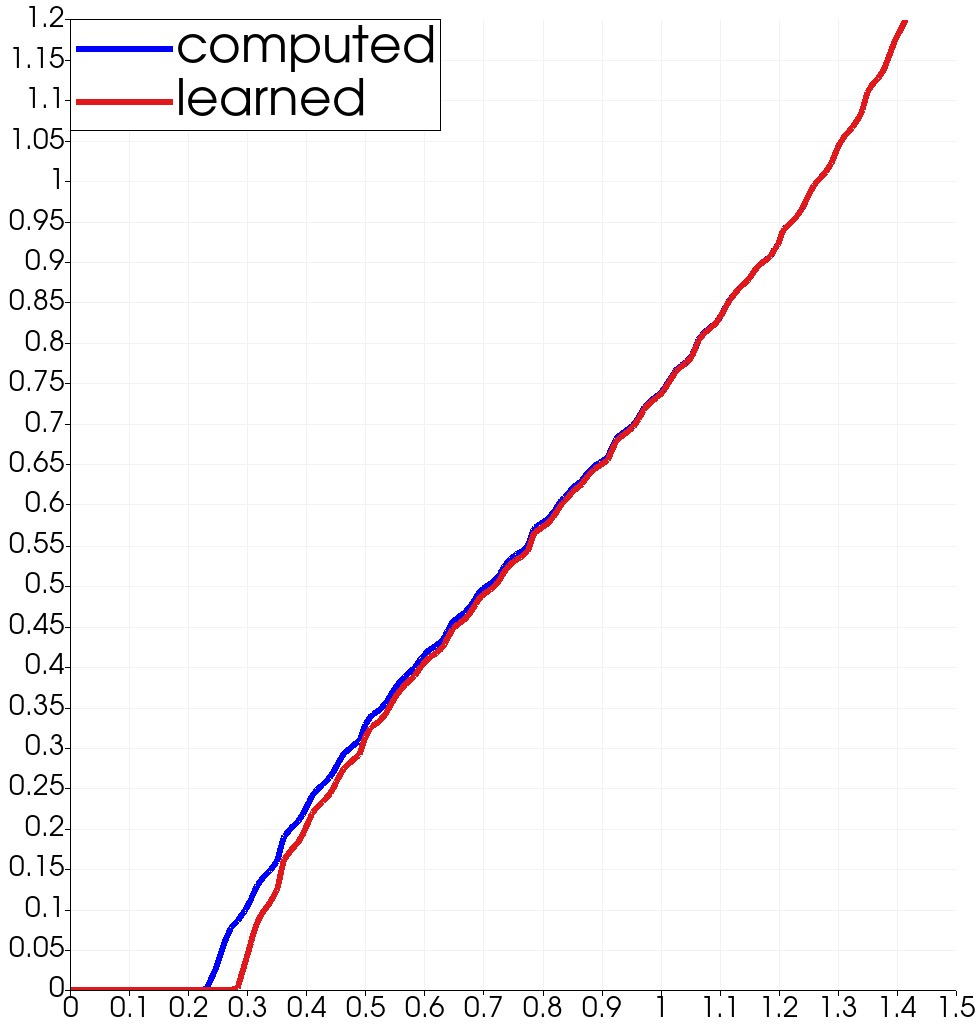}
\end{subfigure}
\begin{subfigure}{0.49\textwidth}
	\centering
	\includegraphics[width=.7\linewidth]{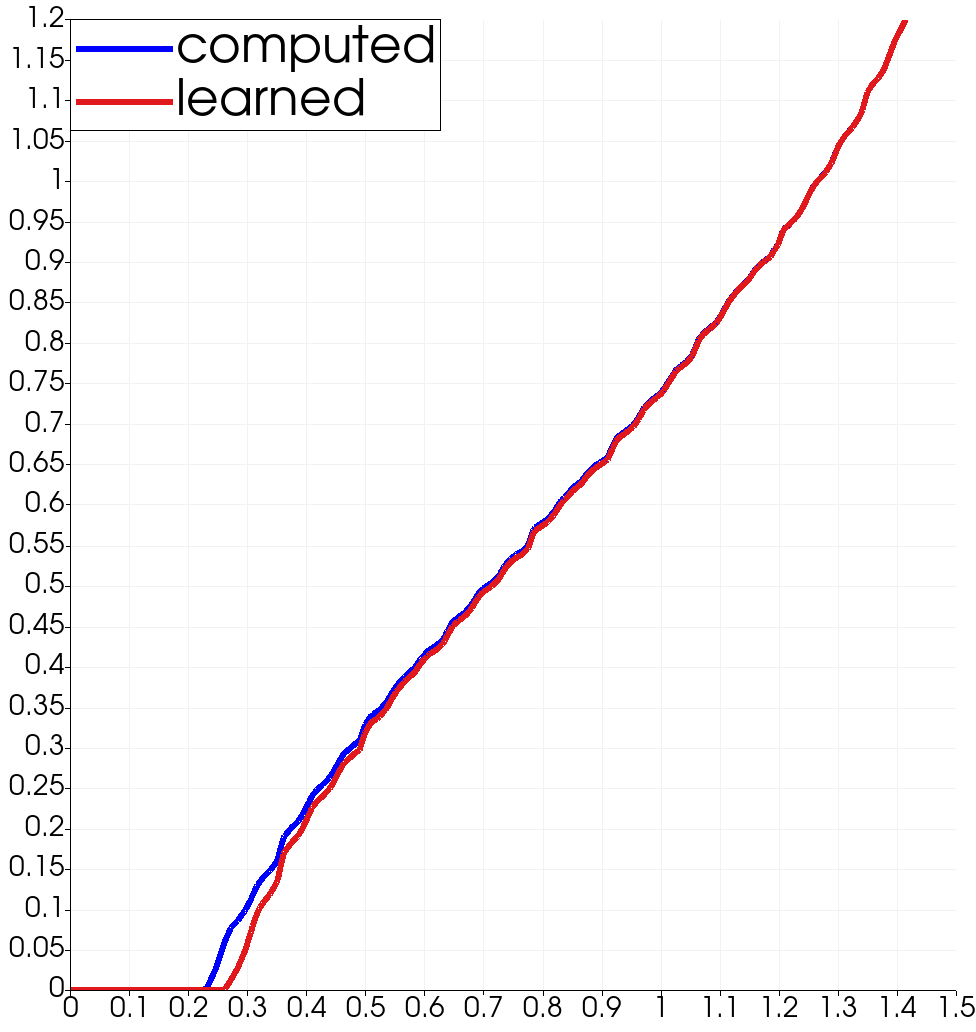}
\end{subfigure}
\caption{Original (blue) and surrogate (red) reconstructed solution along the line $y = x$ for $n_s = 2^4$ (top left), $n_s = 3^4$ (top right), $n_s = 4^4$ (bottom left) and $n_s = 5^4$ (bottom right)} \label{fig_pme_2d_overline}
\end{figure}

The solution obtained using the surrogate models trained with $n_s = 2^4, 3^4, 4^4$ and $5^4$ are reported on Figure 
\ref{fig_pme_2d_surrogate}, together with solution profiles  along the line $x=y$ shown on Figure \ref{fig_pme_2d_overline}. 

\begin{table}[h]
\begin{center}
\begin{tabular}{|c|c|c|c|c|}
\hline
$n_s$ & $2^4$ & $3^4$ & $4^4$ & $5^4$ \\\hline
Rel. error & $33\%$   &  $8\%$     &   $4\%$   &  $3\%$   \\\hline
\end{tabular}
\caption{Relative $L^2(\O)$ error between surrogate model and true models.}
\label{tab:err2d}
\end{center}
\end{table}
We observe that the error decreases as the number of sampling points grows, with the relative $L^2(\O)$ being below $10\%$ starting from $n_s = 3^4$.

\begin{figure}[h]
\centering
\begin{subfigure}{0.49\textwidth}
	\centering
	\includegraphics[width=.7\linewidth]{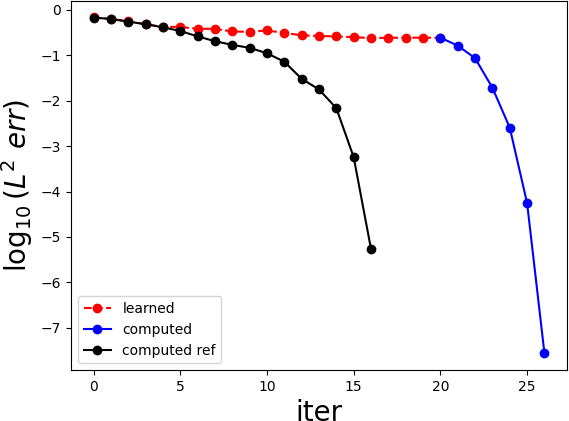}
\end{subfigure}
\begin{subfigure}{0.49\textwidth}
	\centering
	\includegraphics[width=.7\linewidth]{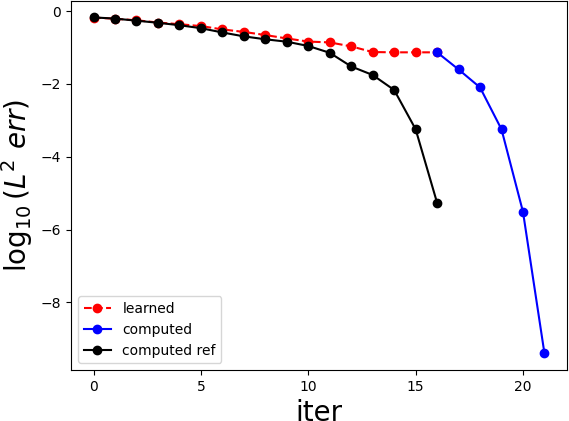}
\end{subfigure}
\begin{subfigure}{0.49\textwidth}
	\centering
	\includegraphics[width=.7\linewidth]{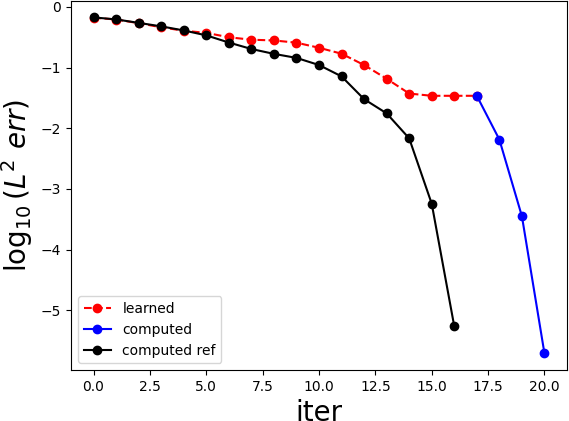}
\end{subfigure}
\caption{$L^2$ error between the reference solution and the current iterate of Newton's method for \eqref{eq:approxsubstlearn} using zero initial guess (red), and \eqref{eq:nonlincoarsesubstruct} with either initial guess set to zero (black) or provided by the approximate solution of \eqref{eq:approxsubstlearn}. Top row from left to right, $n_s = 2^4, 3^4$, bottom row $n_s =4^4$} \label{newton_init}
\end{figure}

\kb{Let us recall that the evaluation of the  residual function \eqref{eq:FH} based on the original coarse
$\dtn_{H,i}$ operator requires the solution of the family of local nonlinear finite element systems, which is computationally demanding. In contrast, the inference time of $\wdtn_{H,i}$ is very low, allowing for fast computation of the  solution to the surrogate system
\mb{\eqref{eq:approxsubstlearn}}. \red{(X)}. Although the surrogate solution $\widetilde{g}_H$ may not be sufficiently accurate, it can serve to initialize the nonlinear solver for the original coarse problem \eqref{eq:nonlincoarsesubstruct}.
This idea} is illustrated by Figure \ref{newton_init}. Figure \ref{newton_init} reports convergence of the relative $L^2(\O)$ error between a very accurate solution to \eqref{eq:nonlincoarsesubstruct} 
and the current Newton's method iterate obtained using either using $\wdtn_{H,i}$ or $\dtn_{H,i}$ with different initialization. 
\kb{The black curve shows convergence of Newton's method for the original problem \eqref{eq:nonlincoarsesubstruct} starting from zero initial guess and is the same for all three sub-figures. Similarly, the red curves depict convergence of Newton's method for the surrogate problem \eqref{eq:approxsubstlearn}  for various number $n_s = 2^4, 3^4$ and $4^4$, again starting for zero initial guess. The stagnation of the red curve is results the approximation error associated to $\wdtn_{H,i}$. Finally, the blue curve shows convergence of Newton's method for \eqref{eq:nonlincoarsesubstruct} starting from the surrogate solution. For the original coarse problem, convergence up to the tolerance of $10^{-5}$ is obtained within 17 iterations.
In contrast, the Newton's method initialized by the surrogate solution requires up to 4 times fewer iterations in order to achieve the similar accuracy. In particular, there are only  7, 5 or 4 iterations needed for  $n_s = 2^4, 3^4$ and $4^4$, respectively. Note that even when using a very crude $n_s = 2$ surrogate model, a significant reduction of Newton's step is achieved.}

\section{Conclusion and outlook}\label{sec:conclusion}
We have presented the framework that combines the Multi-scale Finite Element (MsFEM) method with techniques from Machine Learning, providing an extension of the former to the case of nonlinear PDEs with highly oscilatory coefficients. Our approach relies on learning Dirichlet-to-Neumann (DtN) maps acting between some low-dimensional (coarse) spaces. The surrogate DtN operators are further combined together within a global substructured formulation solved by Newton's method. The optimization process involves fitting values and derivatives of the true coarse DtN, where incorporating derivative information into the loss function plays a pivotal role in accuracy of the learned operator. In addition, we weakly enforce some monotonicity properties of the original model, which improves the performance of Newton's method. 

Numerical experiments, performed in 1D and 2D and involving challenging $p-$Laplace and degenerate diffusion equations, have shown promising results. With just a few training points by dimension, the substitution model can approach the solution with an accuracy of a few percent. A further improvement in the accuracy can be achieved by using the ``learned" solution as an initial guess for Newton's method applied to the original (not learned) substructured formulation, allowing us to enter directly into the region of quadratic convergence and obtain a more accurate ``learned" solution. 
Although the our numerical experiment has been so far limited to the case of periodically distributed coefficients, the extension to the non-periodic case is trivial and will be carried out in the future. 
Future work also involves
incorporating techniques from Principal Component Analysis (PCA) to reduce complexity associated to the derivatives of the DtN. Additionally, the inclusion of higher-order derivatives could increase model accuracy. 
Higher-order optimizers like BFGS, which incorporates  curvature/Hessian information into the loss function, may improve the convergence of the loss function.

\bibliography{bibliography}
\bibliographystyle{spmpsci}

\end{document}